\documentclass[10pt,a4,twoside,rm]{article}
\usepackage[left=4.5cm,top=4cm,right=4.5cm,bottom=2cm]{geometry}
\usepackage{amsmath, amsthm, amssymb}
\usepackage{graphicx}
\usepackage[all, knot]{xy}
\usepackage{amssymb}
\usepackage{fancyhdr}
\usepackage{titlesec}
\usepackage{youngtab}
\usepackage{color}
\usepackage{hyperref}

\def\linew(#1){%
/xywidth{#1 setlinewidth}def}

\def\action{\xy 0;/r1pc/:
(-12,0); (12,0) **\dir{-};
(-7,1); (-7,-1); **\dir{.};
(-2,1); (-2,-1) **\dir{.};
(3,1); (3,-1) **\dir{.};
(8,1); (8,-1) **\dir{.};
(0,.5); (0,-.5) **\dir{-};
(6,.5); (6,-.5) **\dir{-};
(10,.5); (10,-.5) **\dir{-};
(-4,.5); (-4,-.5) **\dir{-};
(-10,.5); (-10,-.5) **\dir{-};
(13,0)*{\cdots};
(-13,0)*{\cdots};
{\ar@/^1pc/(0.2,1)*{};(5.8,1)*{}};
{\ar@/^1pc/(6.2,1)*{};(9.8,1)*{}};
{\ar@/^1pc/(-3.8,1)*{};(-0.2,1)*{}};
{\ar@/^1pc/(-9.8,1)*{};(-4.2,1)*{}};
(0,-1.5)*{\scriptstyle \bf 0};
(-4,-1.5)*{\scriptstyle \bf -4};
(-10,-1.5)*{\scriptstyle \bf-10};
(6,-1.5)*{ \scriptstyle \bf 6};
(10,-1.5)*{\scriptstyle \bf 10};
(-7,-1.5)*{ \scriptstyle -7};
(-2,-1.5)*{\scriptstyle -2};
(3,-1.5)*{ \scriptstyle 3};
(8,-1.5)*{\scriptstyle  8};
\endxy}

\def\indexation{\xy 0;/r.7pc/:
(-20,0); (16,0) **\dir{-};
(-17,1); (-17,-1); **\dir{.};
(-12,1); (-12,-1); **\dir{.};
(-7,1); (-7,-1); **\dir{.};
(-2,1); (-2,-1) **\dir{.};
(3,1); (3,-1) **\dir{.};
(8,1); (8,-1) **\dir{.};
(13,1); (13,-1) **\dir{.};
(1,.5); (1,-.5) **\dir{-};
(5,.5); (5,-.5) **\dir{-};
(11,.5); (11,-.5) **\dir{-};
(15,.5); (15,-.5) **\dir{-};
(-5,.5); (-5,-.5) **\dir{-};
(-9,.5); (-9,-.5) **\dir{-};
(-15,.5); (-15,-.5) **\dir{-};
(-19,.5); (-19,-.5) **\dir{-};
(1,1.5)*{\scriptscriptstyle \bf 1};
(5,1.5)*{\scriptscriptstyle \bf 5};
(11,1.5)*{\scriptscriptstyle \bf 11};
(15,1.5)*{\scriptscriptstyle \bf 15};
(-5,1.5)*{\scriptscriptstyle \bf -5};
(-9,1.5)*{\scriptscriptstyle \bf -9};
(-15,1.5)*{\scriptscriptstyle \bf -15};
(-19,1.5)*{\scriptscriptstyle \bf -19};
(-19,-1.5)*{\scriptstyle \bs{\lambda}_1};
(15,-1.5)*{\scriptstyle\bs{\lambda}_2};
(-15,-1.5)*{\scriptstyle\bs{\lambda}_3};
(11,-1.5)*{\scriptstyle\bs{\lambda}_4};
(-9,-1.5)*{\scriptstyle\bs{\lambda}_5};
(5,-1.5)*{\scriptstyle\bs{\lambda}_6};
(-5,-1.5)*{\scriptstyle\bs{\lambda}_7};
(1,-1.5)*{\scriptstyle\bs{\lambda}_8};
(-7,1.5)*{ \scriptscriptstyle -7};
(-2,1.5)*{\scriptscriptstyle -2};
(3,1.5)*{ \scriptscriptstyle 3};
(8,1.5)*{\scriptscriptstyle  8};
(13,1.5)*{\scriptscriptstyle  13};
(-12,1.5)*{\scriptscriptstyle -12};
(-17,1.5)*{\scriptscriptstyle -17};
\endxy}

\def\eje{\xy 0;/r.8pc/:
(0.5,0.5)*{ \sub{s} = ( };
(2,0); (2,1) **\dir{-};
(2,0); (7,0) **\dir{-};
(2,1); (7,1) **\dir{-};
(3,0); (3,1) **\dir{-};
(4,0); (4,1) **\dir{-};
(5,0); (5,1) **\dir{-};
(6,0); (6,1) **\dir{-};
(7,0); (7,1) **\dir{-};
(2.5,0.5)*{ 1};
(3.5,0.5)*{ 2};
(4.5,0.5)*{3 };
(5.5,0.5)*{ 4};
(6.5,0.5)*{ 6};
(7.5,0)*{ ,};
(8,0); (8,1) **\dir{-};
(9,0); (9,1) **\dir{-};
(8,1); (9,1) **\dir{-};
(8,0); (9,0) **\dir{-};
(8.5,0.5)*{ 5};
(9.5,0.5)*{ ) };
(15.5,0.5)*{ \sub{t} = ( };
(17,0); (17,1) **\dir{-};
(17,0); (22,0) **\dir{-};
(17,1); (22,1) **\dir{-};
(18,0); (18,1) **\dir{-};
(19,0); (19,1) **\dir{-};
(20,0); (20,1) **\dir{-};
(21,0); (21,1) **\dir{-};
(22,0); (22,1) **\dir{-};
(17.5,0.5)*{ 1};
(18.5,0.5)*{ 2};
(19.5,0.5)*{3 };
(20.5,0.5)*{ 4};
(21.5,0.5)*{ 5};
(22.5,0)*{ ,};
(23,0); (23,1) **\dir{-};
(24,0); (24,1) **\dir{-};
(23,1); (24,1) **\dir{-};
(23,0); (24,0) **\dir{-};
(23.5,0.5)*{ 6};
(24.5,0.5)*{ ) };
\endxy}

\def\dibujo{\xy 0;/r.6pc/:
(-1,-1); (0,0) **\dir{-};
(-1,-1); (0,-2) **\dir{-};
(-1,-3); (0,-2) **\dir{-};
(-1,-3); (0,-4) **\dir{-};
(-2,-6); (0,-4) **\dir{-};
(-2,-6); (3,-11) **\dir{-};
(-8,-22); (3,-11) **\dir{-};
(0,0)*{ \bullet};
(-1,-1)*{ \bullet};
(1,-1)*{ \bullet};
(-2,-2)*{ \bullet};
(2,-2)*{ \bullet};
(0,-2)*{ \bullet};
(-3,-3)*{ \bullet};
(3,-3)*{ \bullet};
(-1,-3)*{ \bullet};
(1,-3)*{ \bullet};
(4,-4)*{ \bullet};
(-4,-4)*{ \bullet};
(2,-4)*{ \bullet};
(-2,-4)*{ \bullet};
(0,-4)*{ \bullet};
(-3,-5)*{ \bullet};
(3,-5)*{ \bullet};
(-1,-5)*{ \bullet};
(1,-5)*{ \bullet};
(-5,-5)*{ \bullet};
(5,-5)*{ \bullet};
(4,-6)*{ \bullet};
(-4,-6)*{ \bullet};
(2,-6)*{ \bullet};
(-2,-6)*{ \bullet};
(0,-6)*{ \bullet};
(6,-6)*{ \bullet};
(-6,-6)*{ \bullet};
(-3,-7)*{ \bullet};
(3,-7)*{ \bullet};
(-1,-7)*{ \bullet};
(1,-7)*{ \bullet};
(-5,-7)*{ \bullet};
(5,-7)*{ \bullet};
(-7,-7)*{ \bullet};
(7,-7)*{ \bullet};
(4,-8)*{ \bullet};
(-4,-8)*{ \bullet};
(2,-8)*{ \bullet};
(-2,-8)*{ \bullet};
(0,-8)*{ \bullet};
(6,-8)*{ \bullet};
(-6,-8)*{ \bullet};
(8,-8)*{ \bullet};
(-8,-8)*{ \bullet};
(-3,-9)*{ \bullet};
(3,-9)*{ \bullet};
(-1,-9)*{ \bullet};
(1,-9)*{ \bullet};
(-5,-9)*{ \bullet};
(5,-9)*{ \bullet};
(-7,-9)*{ \bullet};
(7,-9)*{ \bullet};
(-9,-9)*{ \bullet};
(9,-9)*{ \bullet};
(4,-10)*{ \bullet};
(-4,-10)*{ \bullet};
(2,-10)*{ \bullet};
(-2,-10)*{ \bullet};
(0,-10)*{ \bullet};
(6,-10)*{ \bullet};
(-6,-10)*{ \bullet};
(8,-10)*{ \bullet};
(-8,-10)*{ \bullet};
(10,-10)*{ \bullet};
(-10,-10)*{ \bullet};
(-3,-11)*{ \bullet};
(3,-11)*{ \bullet};
(-1,-11)*{ \bullet};
(1,-11)*{ \bullet};
(-5,-11)*{ \bullet};
(5,-11)*{ \bullet};
(-7,-11)*{ \bullet};
(7,-11)*{ \bullet};
(-9,-11)*{ \bullet};
(9,-11)*{ \bullet};
(-11,-11)*{ \bullet};
(11,-11)*{ \bullet};
(4,-12)*{ \bullet};
(-4,-12)*{ \bullet};
(2,-12)*{ \bullet};
(-2,-12)*{ \bullet};
(0,-12)*{ \bullet};
(6,-12)*{ \bullet};
(-6,-12)*{ \bullet};
(8,-12)*{ \bullet};
(-8,-12)*{ \bullet};
(10,-12)*{ \bullet};
(-10,-12)*{ \bullet};
(12,-12)*{ \bullet};
(-12,-12)*{ \bullet};
(-3,-13)*{ \bullet};
(3,-13)*{ \bullet};
(-1,-13)*{ \bullet};
(1,-13)*{ \bullet};
(-5,-13)*{ \bullet};
(5,-13)*{ \bullet};
(-7,-13)*{ \bullet};
(7,-13)*{ \bullet};
(-9,-13)*{ \bullet};
(9,-13)*{ \bullet};
(-11,-13)*{ \bullet};
(11,-13)*{ \bullet};
(-13,-13)*{ \bullet};
(13,-13)*{ \bullet};
(-14,-14)*{ \bullet};
(-12,-14)*{ \bullet};
(-10,-14)*{ \bullet};
(-8,-14)*{ \bullet};
(-6,-14)*{ \bullet};
(-4,-14)*{ \bullet};
(-2,-14)*{ \bullet};
(0,-14)*{ \bullet};
(2,-14)*{ \bullet};
(4,-14)*{ \bullet};
(6,-14)*{ \bullet};
(8,-14)*{ \bullet};
(10,-14)*{ \bullet};
(12,-14)*{ \bullet};
(-14,-14)*{ \bullet};
(-3,-15)*{ \bullet};
(3,-15)*{ \bullet};
(-1,-15)*{ \bullet};
(1,-15)*{ \bullet};
(-5,-15)*{ \bullet};
(5,-15)*{ \bullet};
(-7,-15)*{ \bullet};
(7,-15)*{ \bullet};
(-9,-15)*{ \bullet};
(9,-15)*{ \bullet};
(-11,-15)*{ \bullet};
(11,-15)*{ \bullet};
(-13,-15)*{ \bullet};
(13,-15)*{ \bullet};
(-15,-15)*{ \bullet};
(15,-15)*{ \bullet};
(-14,-16)*{ \bullet};
(-12,-16)*{ \bullet};
(-10,-16)*{ \bullet};
(-8,-16)*{ \bullet};
(-6,-16)*{ \bullet};
(-4,-16)*{ \bullet};
(-2,-16)*{ \bullet};
(0,-16)*{ \bullet};
(2,-16)*{ \bullet};
(4,-16)*{ \bullet};
(6,-16)*{ \bullet};
(8,-16)*{ \bullet};
(10,-16)*{ \bullet};
(12,-16)*{ \bullet};
(-14,-16)*{ \bullet};
(-16,-16)*{ \bullet};
(16,-16)*{ \bullet};
(-3,-17)*{ \bullet};
(3,-17)*{ \bullet};
(-1,-17)*{ \bullet};
(1,-17)*{ \bullet};
(-5,-17)*{ \bullet};
(5,-17)*{ \bullet};
(-7,-17)*{ \bullet};
(7,-17)*{ \bullet};
(-9,-17)*{ \bullet};
(9,-17)*{ \bullet};
(-11,-17)*{ \bullet};
(11,-17)*{ \bullet};
(-13,-17)*{ \bullet};
(13,-17)*{ \bullet};
(-15,-17)*{ \bullet};
(15,-17)*{ \bullet};
(-17,-17)*{ \bullet};
(17,-17)*{ \bullet};
(-14,-18)*{ \bullet};
(-12,-18)*{ \bullet};
(-10,-18)*{ \bullet};
(-8,-18)*{ \bullet};
(-6,-18)*{ \bullet};
(-4,-18)*{ \bullet};
(-2,-18)*{ \bullet};
(0,-18)*{ \bullet};
(2,-18)*{ \bullet};
(4,-18)*{ \bullet};
(6,-18)*{ \bullet};
(8,-18)*{ \bullet};
(10,-18)*{ \bullet};
(12,-18)*{ \bullet};
(-14,-18)*{ \bullet};
(-16,-18)*{ \bullet};
(16,-18)*{ \bullet};
(-18,-18)*{ \bullet};
(18,-18)*{ \bullet};
(-3,-19)*{ \bullet};
(3,-19)*{ \bullet};
(-1,-19)*{ \bullet};
(1,-19)*{ \bullet};
(-5,-19)*{ \bullet};
(5,-19)*{ \bullet};
(-7,-19)*{ \bullet};
(7,-19)*{ \bullet};
(-9,-19)*{ \bullet};
(9,-19)*{ \bullet};
(-11,-19)*{ \bullet};
(11,-19)*{ \bullet};
(-13,-19)*{ \bullet};
(13,-19)*{ \bullet};
(-15,-19)*{ \bullet};
(15,-19)*{ \bullet};
(-17,-19)*{ \bullet};
(17,-19)*{ \bullet};
(-19,-19)*{ \bullet};
(19,-19)*{ \bullet};
(-14,-20)*{ \bullet};
(-12,-20)*{ \bullet};
(-10,-20)*{ \bullet};
(-8,-20)*{ \bullet};
(-6,-20)*{ \bullet};
(-4,-20)*{ \bullet};
(-2,-20)*{ \bullet};
(0,-20)*{ \bullet};
(2,-20)*{ \bullet};
(4,-20)*{ \bullet};
(6,-20)*{ \bullet};
(8,-20)*{ \bullet};
(10,-20)*{ \bullet};
(12,-20)*{ \bullet};
(-14,-20)*{ \bullet};
(-16,-20)*{ \bullet};
(16,-20)*{ \bullet};
(-18,-20)*{ \bullet};
(18,-20)*{ \bullet};
(-20,-20)*{ \bullet};
(20,-20)*{ \bullet};
(14,-14)*{ \bullet};
(14,-16)*{ \bullet};
(14,-18)*{ \bullet};
(14,-20)*{ \bullet};
(-1,-21)*{ \bullet};
(1,-21)*{ \bullet};
(-3,-21)*{ \bullet};
(3,-21)*{ \bullet};
(-1,-21)*{ \bullet};
(1,-21)*{ \bullet};
(-5,-21)*{ \bullet};
(5,-21)*{ \bullet};
(-7,-21)*{ \bullet};
(7,-21)*{ \bullet};
(-9,-21)*{ \bullet};
(9,-21)*{ \bullet};
(-11,-21)*{ \bullet};
(11,-21)*{ \bullet};
(-13,-21)*{ \bullet};
(13,-21)*{ \bullet};
(-15,-21)*{ \bullet};
(15,-21)*{ \bullet};
(-17,-21)*{ \bullet};
(17,-21)*{ \bullet};
(-19,-21)*{ \bullet};
(-21,-21)*{ \bullet};
(21,-21)*{ \bullet};
(19,-21)*{ \bullet};
(-14,-22)*{ \bullet};
(-12,-22)*{ \bullet};
(-10,-22)*{ \bullet};
(-8,-22)*{ \bullet};
(-6,-22)*{ \bullet};
(-4,-22)*{ \bullet};
(-2,-22)*{ \bullet};
(0,-22)*{ \bullet};
(2,-22)*{ \bullet};
(4,-22)*{ \bullet};
(6,-22)*{ \bullet};
(8,-22)*{ \bullet};
(10,-22)*{ \bullet};
(12,-22)*{ \bullet};
(-14,-22)*{ \bullet};
(-16,-22)*{ \bullet};
(16,-22)*{ \bullet};
(-18,-22)*{ \bullet};
(18,-22)*{ \bullet};
(-20,-22)*{ \bullet};
(20,-22)*{ \bullet};
(-22,-22)*{ \bullet};
(22,-22)*{ \bullet};
(14,-14)*{ \bullet};
(14,-16)*{ \bullet};
(14,-18)*{ \bullet};
(14,-20)*{ \bullet};
(14,-22)*{ \bullet};
(-12,0); (-12,-22.5) **\dir{.};
(-7,0); (-7,-22.5) **\dir{.};
(-2,0); (-2,-22.5) **\dir{.};
(3,0); (3,-22.5) **\dir{.};
(8,0); (8,-22.5) **\dir{.};
(13,0); (13,-22.5) **\dir{.};
(18,0); (18,-22.5) **\dir{.};
(-17,0); (-17,-22.5) **\dir{.};
\endxy}

\def\dibujos{\xy 0;/r.5pc/:
\xy 0;/r.6pc/:
(-3,-3); (0,0) **\dir{-};
(-3,-3); (3,-9) **\dir{-};
(-1,-13); (3,-9) **\dir{-};
(0,0)*{ \bullet};
(-1,-1)*{ \bullet};
(1,-1)*{ \bullet};
(-2,-2)*{ \bullet};
(2,-2)*{ \bullet};
(0,-2)*{ \bullet};
(-3,-3)*{ \bullet};
(3,-3)*{ \bullet};
(-1,-3)*{ \bullet};
(1,-3)*{ \bullet};
(4,-4)*{ \bullet};
(-4,-4)*{ \bullet};
(2,-4)*{ \bullet};
(-2,-4)*{ \bullet};
(0,-4)*{ \bullet};
(-3,-5)*{ \bullet};
(3,-5)*{ \bullet};
(-1,-5)*{ \bullet};
(1,-5)*{ \bullet};
(-5,-5)*{ \bullet};
(5,-5)*{ \bullet};
(4,-6)*{ \bullet};
(-4,-6)*{ \bullet};
(2,-6)*{ \bullet};
(-2,-6)*{ \bullet};
(0,-6)*{ \bullet};
(6,-6)*{ \bullet};
(-6,-6)*{ \bullet};
(-3,-7)*{ \bullet};
(3,-7)*{ \bullet};
(-1,-7)*{ \bullet};
(1,-7)*{ \bullet};
(-5,-7)*{ \bullet};
(5,-7)*{ \bullet};
(-7,-7)*{ \bullet};
(7,-7)*{ \bullet};
(4,-8)*{ \bullet};
(-4,-8)*{ \bullet};
(2,-8)*{ \bullet};
(-2,-8)*{ \bullet};
(0,-8)*{ \bullet};
(6,-8)*{ \bullet};
(-6,-8)*{ \bullet};
(8,-8)*{ \bullet};
(-8,-8)*{ \bullet};
(-3,-9)*{ \bullet};
(3,-9)*{ \bullet};
(-1,-9)*{ \bullet};
(1,-9)*{ \bullet};
(-5,-9)*{ \bullet};
(5,-9)*{ \bullet};
(-7,-9)*{ \bullet};
(7,-9)*{ \bullet};
(-9,-9)*{ \bullet};
(9,-9)*{ \bullet};
(4,-10)*{ \bullet};
(-4,-10)*{ \bullet};
(2,-10)*{ \bullet};
(-2,-10)*{ \bullet};
(0,-10)*{ \bullet};
(6,-10)*{ \bullet};
(-6,-10)*{ \bullet};
(8,-10)*{ \bullet};
(-8,-10)*{ \bullet};
(10,-10)*{ \bullet};
(-10,-10)*{ \bullet};
(-3,-11)*{ \bullet};
(3,-11)*{ \bullet};
(-1,-11)*{ \bullet};
(1,-11)*{ \bullet};
(-5,-11)*{ \bullet};
(5,-11)*{ \bullet};
(-7,-11)*{ \bullet};
(7,-11)*{ \bullet};
(-9,-11)*{ \bullet};
(9,-11)*{ \bullet};
(-11,-11)*{ \bullet};
(11,-11)*{ \bullet};
(4,-12)*{ \bullet};
(-4,-12)*{ \bullet};
(2,-12)*{ \bullet};
(-2,-12)*{ \bullet};
(0,-12)*{ \bullet};
(6,-12)*{ \bullet};
(-6,-12)*{ \bullet};
(8,-12)*{ \bullet};
(-8,-12)*{ \bullet};
(10,-12)*{ \bullet};
(-10,-12)*{ \bullet};
(12,-12)*{ \bullet};
(-12,-12)*{ \bullet};
(-3,-13)*{ \bullet};
(3,-13)*{ \bullet};
(-1,-13)*{ \bullet};
(1,-13)*{ \bullet};
(-5,-13)*{ \bullet};
(5,-13)*{ \bullet};
(-7,-13)*{ \bullet};
(7,-13)*{ \bullet};
(-9,-13)*{ \bullet};
(9,-13)*{ \bullet};
(-11,-13)*{ \bullet};
(11,-13)*{ \bullet};
(-13,-13)*{ \bullet};
(13,-13)*{ \bullet};
(-12,0); (-12,-14.5) **\dir{.};
(-7,0); (-7,-14.5) **\dir{.};
(-2,0); (-2,-14.5) **\dir{.};
(3,0); (3,-14.5) **\dir{.};
(8,0); (8,-14.5) **\dir{.};
(13,0); (13,-14.5) **\dir{.};
{\ar@{->} (-14,-15)*{};(14,-15)*{}};
(-12,-15.5)*{\scriptscriptstyle -12};
(-7,-15.5)*{\scriptscriptstyle -7};
(-2,-15.5)*{\scriptscriptstyle -2};
(3,-15.5)*{\scriptscriptstyle 3};
(8,-15.5)*{\scriptscriptstyle 8};
(13,-15.5)*{\scriptscriptstyle 13};
(0,-14.5)*{\scriptscriptstyle  \text{weight}};
\endxy}

\def\dibujoss{\xy 0;/r.6pc/:
(0,0); (-13,-13) **\dir{-};
(-11,-13); (-12,-12) **\dir{-};
(-7,-7); (-2,-12) **\dir{-};
(-3,-13); (-2,-12) **\dir{-};
(-1,-13); (-2,-12) **\dir{-};
(-2,-2); (9,-13) **\dir{-};
(3,-7); (-2,-12) **\dir{-};
(7,-13); (8,-12) **\dir{-};
(0,0)*{ \bullet};
(-1,-1)*{ \bullet};
(1,-1)*{ \bullet};
(-2,-2)*{ \bullet};
(2,-2)*{ \bullet};
(0,-2)*{ \bullet};
(-3,-3)*{ \bullet};
(3,-3)*{ \bullet};
(-1,-3)*{ \bullet};
(1,-3)*{ \bullet};
(4,-4)*{ \bullet};
(-4,-4)*{ \bullet};
(2,-4)*{ \bullet};
(-2,-4)*{ \bullet};
(0,-4)*{ \bullet};
(-3,-5)*{ \bullet};
(3,-5)*{ \bullet};
(-1,-5)*{ \bullet};
(1,-5)*{ \bullet};
(-5,-5)*{ \bullet};
(5,-5)*{ \bullet};
(4,-6)*{ \bullet};
(-4,-6)*{ \bullet};
(2,-6)*{ \bullet};
(-2,-6)*{ \bullet};
(0,-6)*{ \bullet};
(6,-6)*{ \bullet};
(-6,-6)*{ \bullet};
(-3,-7)*{ \bullet};
(3,-7)*{ \bullet};
(-1,-7)*{ \bullet};
(1,-7)*{ \bullet};
(-5,-7)*{ \bullet};
(5,-7)*{ \bullet};
(-7,-7)*{ \bullet};
(7,-7)*{ \bullet};
(4,-8)*{ \bullet};
(-4,-8)*{ \bullet};
(2,-8)*{ \bullet};
(-2,-8)*{ \bullet};
(0,-8)*{ \bullet};
(6,-8)*{ \bullet};
(-6,-8)*{ \bullet};
(8,-8)*{ \bullet};
(-8,-8)*{ \bullet};
(-3,-9)*{ \bullet};
(3,-9)*{ \bullet};
(-1,-9)*{ \bullet};
(1,-9)*{ \bullet};
(-5,-9)*{ \bullet};
(5,-9)*{ \bullet};
(-7,-9)*{ \bullet};
(7,-9)*{ \bullet};
(-9,-9)*{ \bullet};
(9,-9)*{ \bullet};
(4,-10)*{ \bullet};
(-4,-10)*{ \bullet};
(2,-10)*{ \bullet};
(-2,-10)*{ \bullet};
(0,-10)*{ \bullet};
(6,-10)*{ \bullet};
(-6,-10)*{ \bullet};
(8,-10)*{ \bullet};
(-8,-10)*{ \bullet};
(10,-10)*{ \bullet};
(-10,-10)*{ \bullet};
(-3,-11)*{ \bullet};
(3,-11)*{ \bullet};
(-1,-11)*{ \bullet};
(1,-11)*{ \bullet};
(-5,-11)*{ \bullet};
(5,-11)*{ \bullet};
(-7,-11)*{ \bullet};
(7,-11)*{ \bullet};
(-9,-11)*{ \bullet};
(9,-11)*{ \bullet};
(-11,-11)*{ \bullet};
(11,-11)*{ \bullet};
(4,-12)*{ \bullet};
(-4,-12)*{ \bullet};
(2,-12)*{ \bullet};
(-2,-12)*{ \bullet};
(0,-12)*{ \bullet};
(6,-12)*{ \bullet};
(-6,-12)*{ \bullet};
(8,-12)*{ \bullet};
(-8,-12)*{ \bullet};
(10,-12)*{ \bullet};
(-10,-12)*{ \bullet};
(12,-12)*{ \bullet};
(-12,-12)*{ \bullet};
(-3,-13)*{ \bullet};
(3,-13)*{ \bullet};
(-1,-13)*{ \bullet};
(1,-13)*{ \bullet};
(-5,-13)*{ \bullet};
(5,-13)*{ \bullet};
(-7,-13)*{ \bullet};
(7,-13)*{ \bullet};
(-9,-13)*{ \bullet};
(9,-13)*{ \bullet};
(-11,-13)*{ \bullet};
(11,-13)*{ \bullet};
(-13,-13)*{ \bullet};
(13,-13)*{ \bullet};
(-12,0); (-12,-14.5) **\dir{.};
(-7,0); (-7,-14.5) **\dir{.};
(-2,0); (-2,-14.5) **\dir{.};
(3,0); (3,-14.5) **\dir{.};
(8,0); (8,-14.5) **\dir{.};
(13,0); (13,-14.5) **\dir{.};
{\ar@{->} (-14,-15)*{};(14,-15)*{}};
(-12,-15.5)*{\scriptscriptstyle -12};
(-7,-15.5)*{\scriptscriptstyle -7};
(-2,-15.5)*{\scriptscriptstyle -2};
(3,-15.5)*{\scriptscriptstyle 3};
(8,-15.5)*{\scriptscriptstyle 8};
(13,-15.5)*{\scriptscriptstyle 13};
(0,-14.5)*{\scriptscriptstyle  \text{weight}};
\endxy}

\def\pascal{\xy 0;/r2pc/:
(0,1); (0,-3) **\dir{.};
{\ar@{->} (0,.5)*{};(2,.5)*{}};
(4.3,0.5)*{ \scriptstyle \text{central vertical axis}};
(0,0)*{ \bullet};
(-1,-1)*{ \bullet};
(1,-1)*{ \bullet};
(-2,-2)*{ \bullet};
(2,-2)*{ \bullet};
(0,-2)*{ \bullet};
(-3,-3)*{ \bullet};
(3,-3)*{ \bullet};
(-1,-3)*{ \bullet};
(1,-3)*{ \bullet};
{\ar@{->} (-4,-5)*{};(4,-5)*{}};
(-3,-5.5)*{\scriptstyle -3};
(-2,-5.5)*{\scriptstyle -2};
(-1,-5.5)*{\scriptstyle -1};
(0,-5.5)*{\scriptstyle 0};
(1,-5.5)*{\scriptstyle1};
(2,-5.5)*{ \scriptstyle 2};
(3,-5.5)*{\scriptstyle 3};
(0,-4.5)*{\scriptstyle  \text{weight}};
(4,-5.5)*{\ldots};
(-4,-5.5)*{\ldots};
{\ar@{->} (-6,0)*{};(-6,-5)*{}};
(-5.5,0)*{\scriptstyle 0};
(-5.5,-1)*{\scriptstyle 1};
(-5.5,-2)*{\scriptstyle 2};
(-5.5,-3)*{\scriptstyle 3};
(-7,-1.5)*{\scriptstyle  \text{level}};
\endxy
}

\def\mayor{\xy 0;/r1.3pc/:
(0,0)*{ \bullet};
(-1,-1)*{ \bullet};
(1,-1)*{ \bullet};
(-2,-2)*{ \bullet};
(2,-2)*{ \bullet};
(0,-2)*{ \bullet};
(-3,-3)*{ \bullet};
(3,-3)*{ \bullet};
(-1,-3)*{ \bullet};
(1,-3)*{ \bullet};
(4,-4)*{ \bullet};
(-4,-4)*{ \bullet};
(2,-4)*{ \bullet};
(-2,-4)*{ \bullet};
(0,-4)*{ \bullet};
(-3,-5)*{ \bullet};
(3,-5)*{ \bullet};
(-1,-5)*{ \bullet};
(1,-5)*{ \bullet};
(-5,-5)*{ \bullet};
(5,-5)*{ \bullet};
(4,-6)*{ \bullet};
(-4,-6)*{ \bullet};
(2,-6)*{ \bullet};
(-2,-6)*{ \bullet};
(0,-6)*{ \bullet};
(6,-6)*{ \bullet};
(-6,-6)*{ \bullet};
(0,1); (0,-6) **\dir{.};
(0,0); (-1,-1) **\dir{-};
(0,-2); (-1,-1) **\dir{-};
(0,-2); (-1,-3) **\dir{-};
(0,-4); (-1,-3) **\dir{-};
(0,-4); (1,-5) **\dir{-};
(2,-6); (1,-5) **\dir{-};
{\ar@{->} (-7,-8)*{};(7,-8)*{}};
(-6,-8.5)*{\scriptstyle -6};
(-5,-8.5)*{\scriptstyle -5};
(-4,-8.5)*{\scriptstyle -4};
(-3,-8.5)*{\scriptstyle -3};
(-2,-8.5)*{\scriptstyle -2};
(-1,-8.5)*{\scriptstyle -1};
(0,-8.5)*{\scriptstyle 0};
(1,-8.5)*{\scriptstyle1};
(2,-8.5)*{ \scriptstyle 2};
(3,-8.5)*{\scriptstyle 3};
(6,-8.5)*{\scriptstyle 6};
(5,-8.5)*{\scriptstyle 5};
(4,-8.5)*{\scriptstyle 4};
(0,-7.5)*{\scriptstyle  \text{weight}};
{\ar@{->} (-8,0)*{};(-8,-7)*{}};
(-7.5,0)*{\scriptstyle 0};
(-7.5,-1)*{\scriptstyle 1};
(-7.5,-2)*{\scriptstyle 2};
(-7.5,-3)*{\scriptstyle 3};
(-7.5,-4)*{\scriptstyle 4};
(-7.5,-5)*{\scriptstyle 5};
(-7.5,-6)*{\scriptstyle 6};
(-9,-3)*{\scriptstyle  \text{level}};
\endxy }

\def\walls{\xy 0;/r.8pc/:
(-3,-3); (0,0) **\dir{-};
(-3,-3); (3,-9) **\dir{-};
(-1,-13); (3,-9) **\dir{-};
(0,0)*{ \bullet};
(-1,-1)*{ \bullet};
(1,-1)*{ \bullet};
(-2,-2)*{ \bullet};
(2,-2)*{ \bullet};
(0,-2)*{ \bullet};
(-3,-3)*{ \bullet};
(3,-3)*{ \bullet};
(-1,-3)*{ \bullet};
(1,-3)*{ \bullet};
(4,-4)*{ \bullet};
(-4,-4)*{ \bullet};
(2,-4)*{ \bullet};
(-2,-4)*{ \bullet};
(0,-4)*{ \bullet};
(-3,-5)*{ \bullet};
(3,-5)*{ \bullet};
(-1,-5)*{ \bullet};
(1,-5)*{ \bullet};
(-5,-5)*{ \bullet};
(5,-5)*{ \bullet};
(4,-6)*{ \bullet};
(-4,-6)*{ \bullet};
(2,-6)*{ \bullet};
(-2,-6)*{ \bullet};
(0,-6)*{ \bullet};
(6,-6)*{ \bullet};
(-6,-6)*{ \bullet};
(-3,-7)*{ \bullet};
(3,-7)*{ \bullet};
(-1,-7)*{ \bullet};
(1,-7)*{ \bullet};
(-5,-7)*{ \bullet};
(5,-7)*{ \bullet};
(-7,-7)*{ \bullet};
(7,-7)*{ \bullet};
(4,-8)*{ \bullet};
(-4,-8)*{ \bullet};
(2,-8)*{ \bullet};
(-2,-8)*{ \bullet};
(0,-8)*{ \bullet};
(6,-8)*{ \bullet};
(-6,-8)*{ \bullet};
(8,-8)*{ \bullet};
(-8,-8)*{ \bullet};
(-3,-9)*{ \bullet};
(3,-9)*{ \bullet};
(-1,-9)*{ \bullet};
(1,-9)*{ \bullet};
(-5,-9)*{ \bullet};
(5,-9)*{ \bullet};
(-7,-9)*{ \bullet};
(7,-9)*{ \bullet};
(-9,-9)*{ \bullet};
(9,-9)*{ \bullet};
(4,-10)*{ \bullet};
(-4,-10)*{ \bullet};
(2,-10)*{ \bullet};
(-2,-10)*{ \bullet};
(0,-10)*{ \bullet};
(6,-10)*{ \bullet};
(-6,-10)*{ \bullet};
(8,-10)*{ \bullet};
(-8,-10)*{ \bullet};
(10,-10)*{ \bullet};
(-10,-10)*{ \bullet};
(-3,-11)*{ \bullet};
(3,-11)*{ \bullet};
(-1,-11)*{ \bullet};
(1,-11)*{ \bullet};
(-5,-11)*{ \bullet};
(5,-11)*{ \bullet};
(-7,-11)*{ \bullet};
(7,-11)*{ \bullet};
(-9,-11)*{ \bullet};
(9,-11)*{ \bullet};
(-11,-11)*{ \bullet};
(11,-11)*{ \bullet};
(4,-12)*{ \bullet};
(-4,-12)*{ \bullet};
(2,-12)*{ \bullet};
(-2,-12)*{ \bullet};
(0,-12)*{ \bullet};
(6,-12)*{ \bullet};
(-6,-12)*{ \bullet};
(8,-12)*{ \bullet};
(-8,-12)*{ \bullet};
(10,-12)*{ \bullet};
(-10,-12)*{ \bullet};
(12,-12)*{ \bullet};
(-12,-12)*{ \bullet};
(-3,-13)*{ \bullet};
(3,-13)*{ \bullet};
(-1,-13)*{ \bullet};
(1,-13)*{ \bullet};
(-5,-13)*{ \bullet};
(5,-13)*{ \bullet};
(-7,-13)*{ \bullet};
(7,-13)*{ \bullet};
(-9,-13)*{ \bullet};
(9,-13)*{ \bullet};
(-11,-13)*{ \bullet};
(11,-13)*{ \bullet};
(-13,-13)*{ \bullet};
(13,-13)*{ \bullet};
(-12,0); (-12,-14.5) **\dir{.};
(-7,0); (-7,-14.5) **\dir{.};
(-2,0); (-2,-14.5) **\dir{.};
(3,0); (3,-14.5) **\dir{.};
(8,0); (8,-14.5) **\dir{.};
(13,0); (13,-14.5) **\dir{.};
{\ar@{->} (-14,-15)*{};(14,-15)*{}};
(-12,-15.5)*{\scriptscriptstyle -12};
(-7,-15.5)*{\scriptscriptstyle -7};
(-2,-15.5)*{\scriptscriptstyle -2};
(3,-15.5)*{\scriptscriptstyle 3};
(8,-15.5)*{\scriptscriptstyle 8};
(13,-15.5)*{\scriptscriptstyle 13};
(0,-14.5)*{\scriptscriptstyle  \text{weight}};
\endxy
}

\def\maximales{\xy 0;/r.8pc/:
(0,0); (-13,-13) **\dir{-};
(-11,-13); (-12,-12) **\dir{-};
(-7,-7); (-2,-12) **\dir{-};
(-3,-13); (-2,-12) **\dir{-};
(-1,-13); (-2,-12) **\dir{-};
(-2,-2); (9,-13) **\dir{-};
(3,-7); (-2,-12) **\dir{-};
(7,-13); (8,-12) **\dir{-};
(0,0)*{ \bullet};
(-1,-1)*{ \bullet};
(1,-1)*{ \bullet};
(-2,-2)*{ \bullet};
(2,-2)*{ \bullet};
(0,-2)*{ \bullet};
(-3,-3)*{ \bullet};
(3,-3)*{ \bullet};
(-1,-3)*{ \bullet};
(1,-3)*{ \bullet};
(4,-4)*{ \bullet};
(-4,-4)*{ \bullet};
(2,-4)*{ \bullet};
(-2,-4)*{ \bullet};
(0,-4)*{ \bullet};
(-3,-5)*{ \bullet};
(3,-5)*{ \bullet};
(-1,-5)*{ \bullet};
(1,-5)*{ \bullet};
(-5,-5)*{ \bullet};
(5,-5)*{ \bullet};
(4,-6)*{ \bullet};
(-4,-6)*{ \bullet};
(2,-6)*{ \bullet};
(-2,-6)*{ \bullet};
(0,-6)*{ \bullet};
(6,-6)*{ \bullet};
(-6,-6)*{ \bullet};
(-3,-7)*{ \bullet};
(3,-7)*{ \bullet};
(-1,-7)*{ \bullet};
(1,-7)*{ \bullet};
(-5,-7)*{ \bullet};
(5,-7)*{ \bullet};
(-7,-7)*{ \bullet};
(7,-7)*{ \bullet};
(4,-8)*{ \bullet};
(-4,-8)*{ \bullet};
(2,-8)*{ \bullet};
(-2,-8)*{ \bullet};
(0,-8)*{ \bullet};
(6,-8)*{ \bullet};
(-6,-8)*{ \bullet};
(8,-8)*{ \bullet};
(-8,-8)*{ \bullet};
(-3,-9)*{ \bullet};
(3,-9)*{ \bullet};
(-1,-9)*{ \bullet};
(1,-9)*{ \bullet};
(-5,-9)*{ \bullet};
(5,-9)*{ \bullet};
(-7,-9)*{ \bullet};
(7,-9)*{ \bullet};
(-9,-9)*{ \bullet};
(9,-9)*{ \bullet};
(4,-10)*{ \bullet};
(-4,-10)*{ \bullet};
(2,-10)*{ \bullet};
(-2,-10)*{ \bullet};
(0,-10)*{ \bullet};
(6,-10)*{ \bullet};
(-6,-10)*{ \bullet};
(8,-10)*{ \bullet};
(-8,-10)*{ \bullet};
(10,-10)*{ \bullet};
(-10,-10)*{ \bullet};
(-3,-11)*{ \bullet};
(3,-11)*{ \bullet};
(-1,-11)*{ \bullet};
(1,-11)*{ \bullet};
(-5,-11)*{ \bullet};
(5,-11)*{ \bullet};
(-7,-11)*{ \bullet};
(7,-11)*{ \bullet};
(-9,-11)*{ \bullet};
(9,-11)*{ \bullet};
(-11,-11)*{ \bullet};
(11,-11)*{ \bullet};
(4,-12)*{ \bullet};
(-4,-12)*{ \bullet};
(2,-12)*{ \bullet};
(-2,-12)*{ \bullet};
(0,-12)*{ \bullet};
(6,-12)*{ \bullet};
(-6,-12)*{ \bullet};
(8,-12)*{ \bullet};
(-8,-12)*{ \bullet};
(10,-12)*{ \bullet};
(-10,-12)*{ \bullet};
(12,-12)*{ \bullet};
(-12,-12)*{ \bullet};
(-3,-13)*{ \bullet};
(3,-13)*{ \bullet};
(-1,-13)*{ \bullet};
(1,-13)*{ \bullet};
(-5,-13)*{ \bullet};
(5,-13)*{ \bullet};
(-7,-13)*{ \bullet};
(7,-13)*{ \bullet};
(-9,-13)*{ \bullet};
(9,-13)*{ \bullet};
(-11,-13)*{ \bullet};
(11,-13)*{ \bullet};
(-13,-13)*{ \bullet};
(13,-13)*{ \bullet};
(-12,0); (-12,-14.5) **\dir{.};
(-7,0); (-7,-14.5) **\dir{.};
(-2,0); (-2,-14.5) **\dir{.};
(3,0); (3,-14.5) **\dir{.};
(8,0); (8,-14.5) **\dir{.};
(13,0); (13,-14.5) **\dir{.};
{\ar@{->} (-14,-15)*{};(14,-15)*{}};
(-12,-15.5)*{\scriptscriptstyle -12};
(-7,-15.5)*{\scriptscriptstyle -7};
(-2,-15.5)*{\scriptscriptstyle -2};
(3,-15.5)*{\scriptscriptstyle 3};
(8,-15.5)*{\scriptscriptstyle 8};
(13,-15.5)*{\scriptscriptstyle 13};
(0,-14.5)*{\scriptscriptstyle  \text{weight}};
\endxy
}

\def\maximalesdos{\xy 0;/r.8pc/:
(0,0); (-16,-16) **\dir{-};
(-2,-2); (12,-16) **\dir{-};
(-8,-16); (-12,-12) **\dir{-};
(-7,-7); (2,-16) **\dir{-};
(-3,-13); (-2,-12) **\dir{-};
(-3,-13); (-6,-16) **\dir{-};
(3,-7); (-2,-12) **\dir{-};
(4,-16); (8,-12) **\dir{-};
(0,0)*{ \bullet};
(-1,-1)*{ \bullet};
(1,-1)*{ \bullet};
(-2,-2)*{ \bullet};
(2,-2)*{ \bullet};
(0,-2)*{ \bullet};
(-3,-3)*{ \bullet};
(3,-3)*{ \bullet};
(-1,-3)*{ \bullet};
(1,-3)*{ \bullet};
(4,-4)*{ \bullet};
(-4,-4)*{ \bullet};
(2,-4)*{ \bullet};
(-2,-4)*{ \bullet};
(0,-4)*{ \bullet};
(-3,-5)*{ \bullet};
(3,-5)*{ \bullet};
(-1,-5)*{ \bullet};
(1,-5)*{ \bullet};
(-5,-5)*{ \bullet};
(5,-5)*{ \bullet};
(4,-6)*{ \bullet};
(-4,-6)*{ \bullet};
(2,-6)*{ \bullet};
(-2,-6)*{ \bullet};
(0,-6)*{ \bullet};
(6,-6)*{ \bullet};
(-6,-6)*{ \bullet};
(-3,-7)*{ \bullet};
(3,-7)*{ \bullet};
(-1,-7)*{ \bullet};
(1,-7)*{ \bullet};
(-5,-7)*{ \bullet};
(5,-7)*{ \bullet};
(-7,-7)*{ \bullet};
(7,-7)*{ \bullet};
(4,-8)*{ \bullet};
(-4,-8)*{ \bullet};
(2,-8)*{ \bullet};
(-2,-8)*{ \bullet};
(0,-8)*{ \bullet};
(6,-8)*{ \bullet};
(-6,-8)*{ \bullet};
(8,-8)*{ \bullet};
(-8,-8)*{ \bullet};
(-3,-9)*{ \bullet};
(3,-9)*{ \bullet};
(-1,-9)*{ \bullet};
(1,-9)*{ \bullet};
(-5,-9)*{ \bullet};
(5,-9)*{ \bullet};
(-7,-9)*{ \bullet};
(7,-9)*{ \bullet};
(-9,-9)*{ \bullet};
(9,-9)*{ \bullet};
(4,-10)*{ \bullet};
(-4,-10)*{ \bullet};
(2,-10)*{ \bullet};
(-2,-10)*{ \bullet};
(0,-10)*{ \bullet};
(6,-10)*{ \bullet};
(-6,-10)*{ \bullet};
(8,-10)*{ \bullet};
(-8,-10)*{ \bullet};
(10,-10)*{ \bullet};
(-10,-10)*{ \bullet};
(-3,-11)*{ \bullet};
(3,-11)*{ \bullet};
(-1,-11)*{ \bullet};
(1,-11)*{ \bullet};
(-5,-11)*{ \bullet};
(5,-11)*{ \bullet};
(-7,-11)*{ \bullet};
(7,-11)*{ \bullet};
(-9,-11)*{ \bullet};
(9,-11)*{ \bullet};
(-11,-11)*{ \bullet};
(11,-11)*{ \bullet};
(4,-12)*{ \bullet};
(-4,-12)*{ \bullet};
(2,-12)*{ \bullet};
(-2,-12)*{ \bullet};
(0,-12)*{ \bullet};
(6,-12)*{ \bullet};
(-6,-12)*{ \bullet};
(8,-12)*{ \bullet};
(-8,-12)*{ \bullet};
(10,-12)*{ \bullet};
(-10,-12)*{ \bullet};
(12,-12)*{ \bullet};
(-12,-12)*{ \bullet};
(-3,-13)*{ \bullet};
(3,-13)*{ \bullet};
(-1,-13)*{ \bullet};
(1,-13)*{ \bullet};
(-5,-13)*{ \bullet};
(5,-13)*{ \bullet};
(-7,-13)*{ \bullet};
(7,-13)*{ \bullet};
(-9,-13)*{ \bullet};
(9,-13)*{ \bullet};
(-11,-13)*{ \bullet};
(11,-13)*{ \bullet};
(-13,-13)*{ \bullet};
(13,-13)*{ \bullet};
(4,-14)*{ \bullet};
(-4,-14)*{ \bullet};
(2,-14)*{ \bullet};
(-2,-14)*{ \bullet};
(0,-14)*{ \bullet};
(6,-14)*{ \bullet};
(-6,-14)*{ \bullet};
(8,-14)*{ \bullet};
(-8,-14)*{ \bullet};
(10,-14)*{ \bullet};
(-10,-14)*{ \bullet};
(12,-14)*{ \bullet};
(-12,-14)*{ \bullet};
(14,-14)*{ \bullet};
(-14,-14)*{ \bullet};
(-3,-15)*{ \bullet};
(3,-15)*{ \bullet};
(-1,-15)*{ \bullet};
(1,-15)*{ \bullet};
(-5,-15)*{ \bullet};
(5,-15)*{ \bullet};
(-7,-15)*{ \bullet};
(7,-15)*{ \bullet};
(-9,-15)*{ \bullet};
(9,-15)*{ \bullet};
(-11,-15)*{ \bullet};
(11,-15)*{ \bullet};
(-13,-15)*{ \bullet};
(13,-15)*{ \bullet};
(-15,-15)*{ \bullet};
(15,-15)*{ \bullet};
(4,-16)*{ \bullet};
(-4,-16)*{ \bullet};
(2,-16)*{ \bullet};
(-2,-16)*{ \bullet};
(0,-16)*{ \bullet};
(6,-16)*{ \bullet};
(-6,-16)*{ \bullet};
(8,-16)*{ \bullet};
(-8,-16)*{ \bullet};
(10,-16)*{ \bullet};
(-10,-16)*{ \bullet};
(12,-16)*{ \bullet};
(-12,-16)*{ \bullet};
(14,-16)*{ \bullet};
(-14,-16)*{ \bullet};
(16,-16)*{ \bullet};
(-16,-16)*{ \bullet};
(-12,0); (-12,-17.5) **\dir{.};
(-7,0); (-7,-17.5) **\dir{.};
(-2,0); (-2,-17.5) **\dir{.};
(3,0); (3,-17.5) **\dir{.};
(8,0); (8,-17.5) **\dir{.};
(13,0); (13,-17.5) **\dir{.};
(-18,-17); (18,-17) **\dir{-};
(-12,1)*{\scriptscriptstyle -12};
(-7,1)*{\scriptscriptstyle -7};
(-2,1)*{\scriptscriptstyle -2};
(3,1)*{\scriptscriptstyle 3};
(8,1)*{\scriptscriptstyle 8};
(13,1)*{\scriptscriptstyle 13};
(-16,-17.5); (-16,-16.5) **\dir{-};
(-8,-17.5); (-8,-16.5) **\dir{-};
(-6,-17.5); (-6,-16.5) **\dir{-};
(2,-17.5); (2,-16.5) **\dir{-};
(4,-17.5); (4,-16.5) **\dir{-};
(12,-17.5); (12,-16.5) **\dir{-};
(-16,-18.5)*{\scriptscriptstyle \la_1};
(-8,-18.5)*{\scriptscriptstyle \la_3};
(-6,-18.5)*{\scriptscriptstyle \la_5};
(2,-18.5)*{\scriptscriptstyle \la_6};
(4,-18.5)*{\scriptscriptstyle \la_4};
(12,-18.5)*{\scriptscriptstyle \la_2};
\endxy
}

\def\maximales{\xy 0;/r.8pc/:
(0,0); (-13,-13) **\dir{-};
(-11,-13); (-12,-12) **\dir{-};
(-7,-7); (-2,-12) **\dir{-};
(-3,-13); (-2,-12) **\dir{-};
(-1,-13); (-2,-12) **\dir{-};
(-2,-2); (9,-13) **\dir{-};
(3,-7); (-2,-12) **\dir{-};
(7,-13); (8,-12) **\dir{-};
(0,0)*{ \bullet};
(-1,-1)*{ \bullet};
(1,-1)*{ \bullet};
(-2,-2)*{ \bullet};
(2,-2)*{ \bullet};
(0,-2)*{ \bullet};
(-3,-3)*{ \bullet};
(3,-3)*{ \bullet};
(-1,-3)*{ \bullet};
(1,-3)*{ \bullet};
(4,-4)*{ \bullet};
(-4,-4)*{ \bullet};
(2,-4)*{ \bullet};
(-2,-4)*{ \bullet};
(0,-4)*{ \bullet};
(-3,-5)*{ \bullet};
(3,-5)*{ \bullet};
(-1,-5)*{ \bullet};
(1,-5)*{ \bullet};
(-5,-5)*{ \bullet};
(5,-5)*{ \bullet};
(4,-6)*{ \bullet};
(-4,-6)*{ \bullet};
(2,-6)*{ \bullet};
(-2,-6)*{ \bullet};
(0,-6)*{ \bullet};
(6,-6)*{ \bullet};
(-6,-6)*{ \bullet};
(-3,-7)*{ \bullet};
(3,-7)*{ \bullet};
(-1,-7)*{ \bullet};
(1,-7)*{ \bullet};
(-5,-7)*{ \bullet};
(5,-7)*{ \bullet};
(-7,-7)*{ \bullet};
(7,-7)*{ \bullet};
(4,-8)*{ \bullet};
(-4,-8)*{ \bullet};
(2,-8)*{ \bullet};
(-2,-8)*{ \bullet};
(0,-8)*{ \bullet};
(6,-8)*{ \bullet};
(-6,-8)*{ \bullet};
(8,-8)*{ \bullet};
(-8,-8)*{ \bullet};
(-3,-9)*{ \bullet};
(3,-9)*{ \bullet};
(-1,-9)*{ \bullet};
(1,-9)*{ \bullet};
(-5,-9)*{ \bullet};
(5,-9)*{ \bullet};
(-7,-9)*{ \bullet};
(7,-9)*{ \bullet};
(-9,-9)*{ \bullet};
(9,-9)*{ \bullet};
(4,-10)*{ \bullet};
(-4,-10)*{ \bullet};
(2,-10)*{ \bullet};
(-2,-10)*{ \bullet};
(0,-10)*{ \bullet};
(6,-10)*{ \bullet};
(-6,-10)*{ \bullet};
(8,-10)*{ \bullet};
(-8,-10)*{ \bullet};
(10,-10)*{ \bullet};
(-10,-10)*{ \bullet};
(-3,-11)*{ \bullet};
(3,-11)*{ \bullet};
(-1,-11)*{ \bullet};
(1,-11)*{ \bullet};
(-5,-11)*{ \bullet};
(5,-11)*{ \bullet};
(-7,-11)*{ \bullet};
(7,-11)*{ \bullet};
(-9,-11)*{ \bullet};
(9,-11)*{ \bullet};
(-11,-11)*{ \bullet};
(11,-11)*{ \bullet};
(4,-12)*{ \bullet};
(-4,-12)*{ \bullet};
(2,-12)*{ \bullet};
(-2,-12)*{ \bullet};
(0,-12)*{ \bullet};
(6,-12)*{ \bullet};
(-6,-12)*{ \bullet};
(8,-12)*{ \bullet};
(-8,-12)*{ \bullet};
(10,-12)*{ \bullet};
(-10,-12)*{ \bullet};
(12,-12)*{ \bullet};
(-12,-12)*{ \bullet};
(-3,-13)*{ \bullet};
(3,-13)*{ \bullet};
(-1,-13)*{ \bullet};
(1,-13)*{ \bullet};
(-5,-13)*{ \bullet};
(5,-13)*{ \bullet};
(-7,-13)*{ \bullet};
(7,-13)*{ \bullet};
(-9,-13)*{ \bullet};
(9,-13)*{ \bullet};
(-11,-13)*{ \bullet};
(11,-13)*{ \bullet};
(-13,-13)*{ \bullet};
(13,-13)*{ \bullet};
(-12,0); (-12,-14.5) **\dir{.};
(-7,0); (-7,-14.5) **\dir{.};
(-2,0); (-2,-14.5) **\dir{.};
(3,0); (3,-14.5) **\dir{.};
(8,0); (8,-14.5) **\dir{.};
(13,0); (13,-14.5) **\dir{.};
{\ar@{->} (-14,-15)*{};(14,-15)*{}};
(-12,-15.5)*{\scriptscriptstyle -12};
(-7,-15.5)*{\scriptscriptstyle -7};
(-2,-15.5)*{\scriptscriptstyle -2};
(3,-15.5)*{\scriptscriptstyle 3};
(8,-15.5)*{\scriptscriptstyle 8};
(13,-15.5)*{\scriptscriptstyle 13};
(0,-14.5)*{\scriptscriptstyle  \text{weight}};
\endxy
}

\newcommand{\la}{\boldsymbol{\lambda} }

\newcommand{\Std}{\operatorname{Std}}
\newcommand{\Shape}{\operatorname{Shape}}

\newtheorem{teo}{Theorem}[section]
\newtheorem{lem}[teo]{Lemma}

\newtheorem{defi}[teo]{Definition}
\newtheorem{cor}[teo]{Corollary}
\newtheorem{exa}[teo]{Example}
\newtheorem{rem}[teo]{Remark}
\newcommand{\C}{ \mathbb C}
\newcommand{\Z}{ \mathbb Z}
\newcommand{\blob}{b_n(m)}

\newcommand{\R}{R}

\newenvironment{dem}{\noindent \textit{Proof:} }{\quad \hfill $\square$}

\newcommand{\sub}[1]{  {  \mathfrak{#1}}}

\newcommand{\pacol}[1]{ {\text{Par}_2(#1)}  }
\newcommand{\bi}[1]{ \text{Bip}_{1}(#1) }
\newcommand{\bs}[1]{ \boldsymbol{#1}}
\newcommand{\diagramof}[1]{#1}

\newcommand{\addable}[2]{{\mathcal{A}}_{\sub{t}}^{#2}(#1)  }
\newcommand{\removable}[2]{{\mathcal{R}}_{\sub{t}}^{#2}(#1)  }

\newcommand{\num}[1]{ n_{  \mathfrak{s}}(#1)}

\newcommand{\de}[1]{ \mathfrak{d} ({ \mathfrak{#1}})  }

\numberwithin{equation}{section}

\titleformat{\section}[hang]{\sc}{\thesection.}{0.5cm}{\filcenter}
\titleformat{\subsection}[hang]{\bf}{\thesubsection.}{0.2cm}{\filright}

\begin{document}

\title{\bf \normalsize GRADED DECOMPOSITION NUMBERS FOR THE BLOB ALGEBRA}
\author{\sc  david plaza{\thanks{ Supported in part by Beca Doctorado Nacional 2011-CONICYT and FONDECYT grant 11211129}}}
\date{}   \maketitle
\begin{abstract}
\noindent \textsc{Abstract.} We give a closed formula for  the graded decomposition numbers of the blob algebra over a field of characteristic zero at a root of unity.
\end{abstract}

\section{Introduction}
The purpose of this article is to continue the investigation, initiated in \cite{Steen-David}, on a two-parameter
generalization $ b_n=b_n(q,m)$  of the Temperley-Lieb algebra that was introduced by
P. Martin and H. Saleur in \cite{Mat-Sal}, as a way of introducing periodicity in the physical model defining the Temperley-Lieb algebra. This algebra also generalizes the well known diagram calculus on the Temperley-Lieb algebra, actually, this algebra can also be defined via a basis of marked (blobbed) Temperley-Lieb diagrams. For this reason $b_n$ was called the blob algebra in \cite{Mat-Sal}.

\medskip
In \cite{Steen-David}, Ryom-Hansen and the author prove that $b_{n}$ is a $\Z$-graded algebra and construct a graded cellular basis for this algebra. The $\Z$ grading on $b_n$ comes from  seminal work of Brundan and Kleshchev \cite{brundan-klesc} in which they prove that the cyclotomic Hecke algebras of type
$ G(r,1,n)$ are $\Z$-graded, and the known realization of $b_n$ as a quotient of the cyclotomic Hecke algebra of type $G(2,1,n)$ (See \cite{martin-wood}). The existence of a explicit graded cellular basis for  $b_n$  allows to define graded cell and simple $b_n$-modules, parameterized by one-line bipartitions of total degree $n$. Graded dimensions of graded cell and simple modules are related by the graded decomposition numbers. Since graded dimension of the graded cell modules for $b_n$ is known from \cite{Steen-David}, the problem of finding graded dimensions of the irreducible $b_n$-modules is equivalent to the problem of finding graded decomposition numbers for $b_n$. The main goal in this paper is to find such numbers.

\medskip
In the ungraded setting, the decomposition numbers for $b_n$ were determined in \cite{martin-wood1} and \cite{Steen1} using completely different methods.
Our approach is essentially combinatorial, and therefore different from those used in the ungraded case. In order to obtain the graded decomposition numbers for $b_n$ we study a family of positively graded cellular subalgebras $b_n(m,\la) $ of $b_n$, with $\la$ a one-line bipartition of total degree $n$. Graded decompositions numbers for $b_n$ and $b_n(m,\la)$ are closely related. Then, we reduce the main problem in this paper to calculate graded decomposition numbers for $b_n(m,\la)$. Now, since $b_n(m,\la)$ is a positively graded algebra we can define a filtration induced by the grading on $b_n(m,\la)$ for each cell $b_n(m,\la)$-module. This filtration together a counting argument is sufficient to determine the graded decomposition numbers for $b_n(m,\la)$ (and therefore for $b_n$).

\medskip
Let us indicate the layout of the article. In section $2$, we introduce basic notions about graded representation theory, our main emphasis will be the concept of graded cellular algebra. In the following section we define the blob algebra and we also recall the main result in \cite{Steen-David}. To be more precise, we explain how the blob algebra is a graded cellular algebra. In  section $4$ we study the degree function defined on the set of all one-line standard bitableaux in order to make the blob algebra a graded cellular algebra. We give a graphical interpretation for this degree function in terms of the Pascal triangle. This interpretation is the key to show that the algebras $b_n(m,\la)$ are positively graded. In the last section we obtain our main result, giving a closed formula for the graded decomposition numbers for $b_n$.

\medskip
\emph{Acknowledgements.} This work is part of my Ph. D. thesis, and it is a pleasure to
thank S. Ryom-Hansen for his encouragement in this project and for useful suggestions that helped me to improve earlier drafts.

\section{Graded Representation Theory}

Throughout this section we fix an  integral domain $R$ and a  $\Z$-graded $R$-algebra free of finite rank $A$. Let $M$ be a finite dimensional graded $A$-module and let $M=\bigoplus_{k\in \mathbb{Z}} M_k$ its direct sum decomposition in homogeneous components, we define its graded  dimension as the Laurent polynomial
\begin{equation}  \label{ Gradeddimensiondefinition}
\text{dim}_t(M) := \sum_{k\in \mathbb{Z}}(\text{dim} M_k) t^{k} \in \mathbb{Z}[t,t^{-1}]
\end{equation}
For a simple graded $A$-module $L$ let $[M : L\langle k \rangle ]$ be the multiplicity of the simple module $L\langle k \rangle $ as a graded composition factor of $M$ for $k \in  \Z$, where $L\langle k \rangle $  is the graded $A$-module obtained by shifting the grading on $L$ up by $k$. Then the graded decomposition number is
\begin{equation}  \label{ gdn definition}
[M:L]_t := \sum_{k\in \mathbb{Z}}[M : L\langle k \rangle ] t^{k} \in \mathbb{Z}[t,t^{-1}]
\end{equation}

If $M$ is a graded $A$-module let $\underline{M}$ be the ungraded $A$-module obtained by forgetting the grading on $M$. All algebras studied in this paper are graded cellular algebras. Let us briefly recall the definition  and some properties of these algebras.  Graded cellular algebras were introduced by Hu and Mathas in \cite[Definition 2.1]{hu-mathas}, following and extending ideas of Graham and Lehrer \cite{gra ler}.

\begin{defi} \label{cellular algebra} \rm
A graded cell datum for $A$ is a quadruple $ (\Lambda, T, C,\deg) $, where $(\Lambda, \geq)$ is a weight poset, $T(\lambda)$ is a finite set for $\lambda \in \Lambda $, and
$$   \begin{array}{rclcc}
      C: \coprod_{\lambda \in \Lambda} T(\lambda) \times T(\lambda) & \rightarrow & A  &  & \deg : \coprod_{\lambda \in \Lambda} T(\lambda) \rightarrow \Z  \\
           (\mathfrak{s},\mathfrak{t})   &   \rightarrow      &  c_{\sub{st}}^{\lambda}  &  &
     \end{array}
   $$
are two functions such that $C$ is injective and

\begin{description}
  \item[\rm (a)] For $\lambda \in \Lambda $ and $\sub{s,t} \in T(\lambda)$, $c_{\sub{st}}^{\lambda}$ is a homogeneous element of A of degree $\deg c_{\sub{st}}^{\lambda} = \deg(\sub{s}) + \deg (\sub{t})$.
  \item[\rm (b)] The set $\{ c_{\sub{st}}^{\lambda} \mbox{ } | \mbox{ }  \sub{s,t} \in T(\lambda) \text{ for } \lambda \in \Lambda \}$ is an $R$-basis of $A$.
  \item[\rm (c)] The $\R$-linear map $*:A \rightarrow A$ determined by $(c_{\sub{st}}^{\lambda})^{*}=c_{\sub{st}}^{\lambda}$, for all $\lambda \in \Lambda$ and all $\sub{s} ,\sub{t}  \in T(\lambda)$, is an algebra anti-automorphism of $A$
  \item[\rm (d)] If $\sub{s},\sub{t} \in T(\lambda) $, for some $\lambda \in \Lambda$, and $a\in A$ then there exist scalars $r_{\sub{us}}(a) \in \R$ such that  $$  ac_{\sub{st}}^{\lambda} \equiv \sum_{\sub{u} \in T(\lambda)} r_{\sub{us}}(a)  c_{\sub{ut}}^{\lambda}  \mod A^{\lambda}$$
      where $A^{\lambda}$ is the $\R$-submodule of $A$ spanned by $\{ c_{\sub{ab}}^{\mu} \mbox{  } | \mbox{  } \mu \geq \lambda ; \sub{a},\sub{b} \in T(\mu)  \}$
\end{description}
The set $\{ c_{\sub{st}}^{\lambda} \mbox{ } | \mbox{ }  \sub{s,t} \in T(\lambda) \text{ for } \lambda \in \Lambda \}$ is called a graded cellular basis for $A$. If $A$ has a graded cellular basis we say that $A$ is a graded cellular algebra. Omitting the axiom (a) in the above definition we recover the original definition of (ungraded) cellular algebras given in \cite{gra ler}.
\end{defi}

Fix a graded cell datum $ (\Lambda, T, C,\deg) $ for $A$. For $\lambda \in \Lambda$ the graded cell module $C(\lambda)$ is the $R$-module with basis $\{c_{\sub{s}}^{\lambda} \mbox{  }  |  \mbox{  } \sub{s} \in T(\lambda) \}$ and $A$-action given by
$$  ac_{\sub{s}}^{\lambda} = \sum_{\sub{u} \in T(\lambda)} r_{\sub{us}}(a)   c_{\sub{u}}^{\lambda}  $$
where the scalars $r_{\sub{us}}(a) \in R$ are the same that appear  in Definition \ref{cellular algebra} (d). The cell modules $C(\lambda)$ are equipped with a homogeneous bilinear form $\langle \, \cdot  , \cdot \,\rangle_{{\lambda}} $ of degree zero determined by
\begin{equation} \label{forma bilineal}
 c_{\sub{as}}^{\lambda}c_{\sub{tb}}^{\lambda} \equiv \langle c_{\sub{s}}^{\lambda}, c_{\sub{t}}^{\lambda} \rangle_{\lambda} c_{\sub{ab}}^{\lambda} \mod{A^{\lambda}}
\end{equation}
for all $\sub{a,b,s,t} \in T(\lambda)$. The radical of this form
$$   \text{rad }C({\lambda}) = \{  x\in C({\lambda}) \, | \, \langle x,y \rangle_{\lambda} =0 \text{ for all } y\in C({\lambda}) \}        $$
is a graded $A$-submodule of $C({\lambda})$ so that $D({\lambda}) = C({\lambda}) / \text{rad } C({\lambda})$ is a graded $A$-module. (See \cite[Proposition 2.9]{Mathas} and \cite[Lemma 2.7]{hu-mathas}). \\
Let $\Lambda_0 = \{   \lambda \in \Lambda \mbox{  } | \mbox{  } D({\lambda}) \neq 0 \}$.
The next theorems classify the simple graded $A$-modules over a field and describe the respective graded decomposition numbers.
\begin{teo}(\cite[Theorem 2.10]{hu-mathas})
Suppose that R is a field. Then $$\{  D(\lambda) \langle k \rangle  \mbox{ }| \mbox{ } \lambda \in \Lambda_{0} \text{ and } k\in \mathbb{Z} \}$$  is a complete set of pairwise non-isomorphic simple graded $A$-modules.
\end{teo}
\begin{teo} \label{gdn in a cellular algebra}
Suppose that $\lambda \in \Lambda$ and $\mu \in \Lambda_0$. Then
\begin{description}
\item[ \rm (a)] $[C(\lambda):D(\mu)]_t \in \mathbb{N}[t,t^{-1}]$;
\item[ \rm (b)] $[C(\lambda):D(\mu)]_{t=1} = [\underline{C}(\lambda):\underline{D}(\mu)]$;
\item[ \rm (c)] $[C(\mu):D(\mu)]_t=1$  and $[C(\lambda):D(\mu)]_t\neq 0$ only if $\lambda \geq \mu$.
\item[ \rm (d)] $\dim_{t}C(\lambda)= \sum_{\lambda \geq \mu} [C(\lambda):D(\mu)]_t \dim_{t} D(\mu)$
\end{description}
\end{teo}
\begin{dem}
Claims (a), (b) and (c) are \cite[Lemma 2.13]{hu-mathas}. Part (d) is a direct consequence of the definitions.
\end{dem}

A graded $A$-module $M=\bigoplus_{i} M_i$ is positively graded if $M_i=0$ whenever $i<0$. That is, all of the homogeneous elements
of $M$ have non-negative degree. By Definition \ref{cellular algebra}(a), it is straightforward to check that a graded cellular algebra is positively graded if and only if $\deg (\sub{s}) \geq 0 $ for all $\sub{s} \in T(\lambda)$ and $\lambda \in \Lambda$. Consequently, if $A$ is positively graded then so is each cell module of $A$. Suppose that $A$ is positively graded and that $M=\bigoplus_iM_{i}$ is a finite dimensional graded $A$-module. For each $j\in \Z$ let ${\mathcal{G}}_jM= \bigoplus_{i\geq j} M_i$. Since $A$ is positively graded,  ${\mathcal{G}}_jM$ is a graded $A$-submodule of $M$. Let $a$ be maximal and $z$ be minimal such that ${\mathcal{G}}_aM=M$ and ${\mathcal{G}}_zM=0$, respectively. Then the grading filtration of $M$ is the filtration
\begin{equation}
0={\mathcal{G}}_zM \leq  {\mathcal{G}}_{z-1}M \leq \ldots {\mathcal{G}}_aM =M
\end{equation}

We finish this section by relating the graded decomposition numbers of $A$ and some graded subalgebras of $A$. Let $ \mathfrak{e} \in A$ be a homogeneous idempotent and let $A_{\mathfrak{e}}$ denote the subalgebra $\mathfrak{e} A \mathfrak{e}$ of $A$. Then $A_{\mathfrak{e}} $ is a graded subalgebra of $A$ and the inclusion $i : A_{\mathfrak{e}} \hookrightarrow A$ is a homogeneous map of degree zero. We write $\mod-A$ (resp. $\mod-A_{\mathfrak{e}}$) for the category of finite dimensional left $A$-modules (resp. $A_{\mathfrak{e}}$-modules).  We define the functor $f: \mod-A \rightarrow \mod-A_{e}$, where for  $V \in \mod-A$, $fV$ is the subspace $\mathfrak{e}V$ of $V$ regarded as $A_{\mathfrak{e}}$-module.

\begin{teo} \label{gdn iguales teorema}
Let  $L$ and $V$  graded $A$-modules. Assume that $L$ is simple and that $\sub{e}L\neq 0$. Then,
\begin{equation}  \label{gdn iguales}
 [V:L]_{t}=[\sub{e} V: \sub{e} L]_{t}
 \end{equation}
where the left (resp. right) side of (\ref{gdn iguales}) correspond to the graded decomposition number for $A$-modules (resp. $A_{\sub{e}}$-modules).
\end{teo}
\begin{dem}
First, we note that a non-zero homogeneous idempotent must have degree zero. Hence, $\mathfrak{e}V$ is a  graded module, where for a homogeneous element $v\in V$ with $\mathfrak{e}v\neq 0$ we have $\deg (v) = \deg(\mathfrak{e} v)$. The same is true for $\sub{e}L$. Therefore, (\ref{gdn iguales}) follows exactly as in the ungraded case \cite[Appendix A1]{don}.
\end{dem}

\section{The blob algebra}

In this section we introduce the main object of our study, the blob algebra. We shall briefly recall the result from \cite{Steen-David}, where it was proved that the blob algebra is a $\mathbb{Z}$-graded algebra and a graded cellular basis was built for this algebra. The blob algebra was introduced by P. Martin and H. Saleur in \cite{Mat-Sal} as a generalization of the Temperley-Lieb algebra. It is usually defined in terms of a basis of blobbed Temperley-Lieb diagrams and their compositions, from which it derives its name. This blob algebra is isomorphic to an algebra defined by a presentation. From now and for the rest of this  paper we fix  an $l$-th root of unity  $q\in \C$, an integer $m$ with $0<m<l$, and a positive integer $n$. For  $k\in \Z$ we define
$$[k]=[k]_q:= q^{k-1} + q^{k-3} + \ldots + q^{-k+1} \in \C$$
the usual Gaussian coefficient.

\begin{defi} \rm  \label{blob}
The blob algebra $b_n=b_n(m) =b_n(q,m)$ is the associative $\C$-algebra on the generators $U_0,U_1,...,U_{n-1}$ subject to the relations
\begin{align*}
U_{i}^{2}&=-[2] U_i &  &\mbox{if }  1 \leq i < n \\
U_iU_{j}U_{i}&=U_i &   & \mbox{if }   |i-j| =1 \text{ and } 1\leq i,j<n \\
U_iU_{j}&=U_{j}U_i &   & \mbox{if } |i-j| > 1  \text{ and } 0\leq i,j<n \\
U_1U_0U_{1}&=[m+1] U_1 &   &  \\
U_0^{2}&=-[m]U_0 & &
\end{align*}
\end{defi}

We now describe the $\Z$-grading on $b_n(m)$ introduced in \cite{Steen-David}. Assume that
\begin{align}
 \label{qym} q^{4}\neq 1,&  & q^{m} \neq q^{-m}, &  & q^{m} \neq q^{2} q^{-m}, & & q^{-m} \neq q^{2} q^{m}  &
\end{align}
Under this restrictions the algebra $b_n(m)$ is known to be a quotient of the Hecke algebra of type $B$, $H_n(q, Q)$ \cite[Proposition 4.4]{martin-wood}, indeed it is also sometimes referred to as the Temperley-Lieb algebra of type $B$. On the other hand, in \cite{brundan-klesc} J. Brundan and A. Kleshchev constructs isomorphisms between cyclotomic Hecke algebras and
Khovanov-Lauda-Rouquier (KLR) algebras (of type $A$). Since the KLR
algebras are $\Z$-graded, the various cyclotomic Hecke algebras become $\Z$-graded in this way as
well. But $H_n(q, Q)$ is also the cyclotomic Hecke algebra of type
$G(2, 1, n) $ and the basic idea in \cite{Steen-David}, in order to prove that $b_n(m)$ is a $\Z$-graded algebra, was to exploit the result from \cite{brundan-klesc} on this quotient construction.

\medskip
Recall that the quantum characteristic of an element $ q $ of a field $ F$ is the smallest positive integer $k$
such that $1 + q + \ldots +q^{k-1}=0$, setting $k=0$ if no such integer exists.  With our choice of $ q \in \C$
the quantum characteristic of $q$ is $ l$ and we shall always assume that it is odd.
We set $I=\mathbb{Z} / l\mathbb{Z}$ and refer to $ I^n$ as the residue sequences of length $ n$.

\begin{teo} \label{definkl}
Let $k\in \Z$ be such that $2k \equiv m \mod l $. The algebra $b_n(m)$ is isomorphic to the algebra with generators
$$ \{ \psi_1,\cdots , \psi_{n-1} \} \cup  \{ y_1,\cdots , y_{n} \}    \cup
  \{ e(\textbf{i}) \mbox{         }  | \mbox{         } \textbf{i} \in I^n \}                $$
subject to the following relations for $\textbf{i},\textbf{j}\in I^n$  and all admissible $r, s$:
\begin{align}
\label{klblob}  e(\bs{i}) &=0,   &\text{if }  i_1=k    \text{ and }   i_2=k-1    \\
\label{klblob1} e(\bs{i}) &=0,   &\text{if }  i_1=-k  \text{ and }    i_2=-k-1      \\
\label{kl0} e(\bs{i}) &=0,   & \text{ if } i_1\neq k, -k \\
 \label{kl1}y_1e(\bs{i})&  = 0,   & \text{ if }  i_{1}=k,-k \\
\label{kl2}e(\textbf{i})e(\textbf{j})& =\delta_{\textbf{i,j}} e(\textbf{i}), & \\
\label{kl3}\sum_{\textbf{i} \in I^n} e(\textbf{i})& =1, \\
\label{kl4}y_{r}e(\textbf{i})& =e(\textbf{i})y_r ,& \\
\label{kl5}\psi_r e(\textbf{i})& =e(s_r \textbf{i}) \psi_r, \\
\label{kl6}y_ry_s& = y_sy_r,&   &   \\
\label{kl7}\psi_ry_s& = y_s\psi_r,&
 \mbox{ if } s\neq r,r+1 \\
\label{kl8}\psi_r\psi_s& = \psi_s\psi_r,&
 \mbox{ if } |s-r|>1 \\
\label{kl9}\psi_ry_{r+1}e(\textbf{i})& =\left\{ \begin{array}{l}
(y_r\psi_r +1)e(\textbf{i}) \\
y_r\psi_r e(\textbf{i})   \\
\end{array}
\right.  &\begin{array}{l}
\mbox{if   }  i_r=i_{r+1}   \\
\mbox{if   }   i_r \neq i_{r+1}     \\
\end{array} \\
\label{kl10}y_{r+1}\psi_re(\textbf{i})& =\left\{ \begin{array}{l}
(\psi_ry_r +1)e(\textbf{i}) \\
\psi_ry_r e(\textbf{i})   \\
\end{array}
\right. &  \begin{array}{lcc}
\mbox{if   }  i_r=i_{r+1}  \\
\mbox{if   }   i_r \neq i_{r+1}   \\
\end{array}
\end{align}
\begin{align}
\label{kl11}\psi_r^{2}e(\textbf{i})& =\left\{ \begin{array}{l}
0   \\
e(\textbf{i})  \\
(y_{r+1}-y_{r})e(\textbf{i}) \\
(y_{r}-y_{r+1})e(\textbf{i})  \\
\end{array}
\right. &  \begin{array}{l}
\mbox{if   }  i_r=i_{r+1}  \\
\mbox{if   }   i_r \neq i_{r+1} \pm 1    \\
\mbox{if   }  i_{r+1}= i_r+1   \\
\mbox{if   }  i_{r+1}= i_r-1  \\
\end{array} \\
\label{kl12}\psi_r\psi_{r+1}\psi_re(\textbf{i})& =\left\{ \begin{array}{l}
(\psi_{r+1}\psi_r\psi_{r +1} +1)e(\textbf{i})  \\
(\psi_{r+1}\psi_r\psi_{r +1} -1)e(\textbf{i})    \\
(\psi_{r+1}\psi_r\psi_{r +1} )e(\textbf{i})    \\
\end{array}
\right. &  \begin{array}{l}
\mbox{if   }  i_{r+2}=i_r=i_{r+1}-1   \\
\mbox{if   } i_{r+2}=i_r=i_{r+1}+1      \\
\mbox{otherwise   }        \\
\end{array}
\end{align}
where $ s_r := (r,r+1)$ is the simple transposition acting on $ I^n $ by permutation of the coordinates $ r, r+1$. The conditions
$$   \begin{array}{ccccc}
       deg\mbox{  } e(\textbf{i})=0,  &  &  deg\mbox{  } y_r=2,  &  &
 deg\mbox{  }\psi_se(\textbf{i})=-a_{i_{s},i_{s+1}} \end{array}
$$
for $1 \leq r\leq n$,  $1 \leq s \leq n-1 $ and $\textbf{i} \in I^{n}$
define a unique $\mathbb{Z}$-grading on $ b_n (m) $ with degree function $ \deg$, where the matrix $(a_{ij})_{i,j\in I}$, is given by
\begin{align*}
a_{ij}& =\left\{ \begin{array}{l}
2 \\
0   \\
-1   \\
\end{array}
\right.   & \begin{array}{l}
              \mbox{if } i=j \\
              \mbox{if } i\neq j\pm 1\\
              \mbox{if } i= j\pm 1.
            \end{array}
\end{align*}
\end{teo}

We now recall the main result in \cite{Steen-David}, to provide a graded cellular basis for $\blob$. In order to establish it precisely we need to introduce a bit more of notation. Recall that a partition $\lambda$ of $n$ is a sequence of weakly
decreasing nonnegative integers $\lambda =(\lambda_1, \lambda_2, \ldots )$ such that $|\lambda|=\sum_{i\geq 0} \lambda_i=n$. The diagram of $\lambda$  is the set
$$[\lambda] =\{ (i,j) \in  \mathbb{N}\times \mathbb{N}  \mbox{  } | \mbox{   } 1\leq j \leq \lambda_i \text{ and } i\geq 1 \}.$$
It is useful to think of $ [\lambda] $ as an array of boxes in the plane, with the indices following
matrix conventions. Thus the box with label $ (i, j)$ belongs to the $i$'th row and $ j$'th column.
For $\lambda$  a partition of $n$ we denote by $\lambda'$ the partition of $n$ obtained from $\lambda$ by interchanging
its rows and columns. In this paper we only work with two-column partitions, that is partitions $\lambda$ such that  $\lambda_i \leq 2 $ for all $i\geq 1$. Denote by $\pacol{n}$ to the set of all two-column partitions of $n$. A $\lambda$-tableau is a bijection $\tau : [\lambda] \rightarrow \{ 1,\ldots,n \} $. We say that $\tau $ has shape $\lambda$ and write $\Shape(\tau)=\lambda $. We say that a tableau is called standard if the entries in each row and each column are increasing. The set of all standard $\lambda$-tableau is denoted by $\Std(\lambda)$. A bipartition of $n$ is a pair $\bs{\lambda}=(\lambda^{(1)},\lambda^{(2)})$ of partitions such that
$$n=|\lambda^{(1)}|+|\lambda^{(2)}|.$$
The diagram of a bipartition $\bs{\lambda}$ is the set
$$[\bs{\lambda}] =\{ (i,j,k) \in  \mathbb{N}\times
\mathbb{N} \times \{1,2\} \mbox{  } | \mbox{   } 1\leq j \leq \lambda_i^{(k)} \}.$$
We can think of it as a pair of diagrams of  partitions. In this paper we only work with one-line bipartitions, that is
bipartitions $\bs{\lambda} $ such that $\lambda^{(k)}_i=0 $ for all $i\geq 2$ and $k=1,2$. Thus for $ c=1,2$, the $c$'th component of
$ [\bs{\lambda}] $ is
$ \{ (i,j,k) \in [\bs{\lambda}] \,  | \,  k=c\}$.
Let $\bi{n}$ be the set of all one-line bipartitions of $n$.
For $\bs{\lambda}$ a bipartition,
a $\bs{\lambda}$-bitableau is a bijection $\sub{t} : [\bs{\lambda}] \rightarrow \{ 1,\ldots,n  \}$.
We say that $\sub{t}$ has shape $\bs{\lambda} $ and write $\Shape(\sub{t})=\bs{\lambda}$.
A $\bs{\lambda}$-bitableau $\sub{t}$ is called standard if the entries of $\sub{t}$ increase from left
to right in each component.  The set of all standard $\bs{\lambda}$-bitableaux is denoted by $\Std(\bs{\lambda})$
and the union $ \bigcup_{\bs{\lambda}} \Std(\bs{\lambda})$ over all bipartitions of $ n$  is denoted by $\Std(\bs{n})$.

\medskip
Set $\Lambda_n :=\{ -n,-n+2,\ldots,n-2,n \}$. Then there is a bijection
$$  f:\bi{n} \rightarrow \Lambda_n \text{,   } ((a),(b)) \rightarrow a-b. $$
We refer to the elements in $\Lambda_n $ as weights and abusing the notation we do not distinguish between one line bipartitions and weights. Using the above function we may define a total order on $\bi{n}$ as follows.
\begin{defi}
Suppose $\bs{\lambda},\bs{\mu} \in \bi{n}$. We define $\bs{\lambda} \succeq \bs{\mu}$
if either $|f(\bs{\lambda})|< |f(\bs{\mu})|$, or if
$|f(\bs{\lambda})|=|f(\bs{\mu})|$ and $f(\bs{\lambda}) \leq f(\bs{\mu})$.
\end{defi}
This order can be extended to $\Std(\bs{n})$ by considering each bitableaux as a sequence of bipartitions, that is, we put $\sub{s} \succeq \sub{t} $ if and only if \begin{equation}  \label{orden en los bitableaux}
\Shape(\sub{s}_{|_k}) \succeq \Shape(\sub{t}_{|_k}) \qquad \text{ for all  } 1\leq k \leq n
\end{equation}
where for $\sub{s}\in \Std (n)$,  $\sub{s}_{|_k} \in \Std (k)$ is the bitableaux obtained from  $\sub{s}$ by erasing the boxes in $\sub{s}$ with entries greater than $k$. We remark that the order $\succeq$ is not in general a total order for $\Std(n)$. Let $\sub{t}^{\bs{\lambda}}$ be the unique standard $\bs{\lambda}$-bitableau
such that $\sub{t}^{\bs{\lambda}} \succeq \sub{t}$ for all
$\sub{t} \in \Std(\bs{\lambda})$. For  $\bs{\lambda} =(a,b)$, set $c=\min\{a,b\} $. Then
in $t^{\bs{\lambda}}$ the numbers $1,2,\ldots, n$ are located increasingly along the rows according to the following rules:
\begin{enumerate}
  \item even numbers less than or equal to $2c$ are placed in the first component.
  \item odd numbers less than $2c$ are placed in the second component.
  \item numbers greater than $2c$ are placed in the remaining boxes.
\end{enumerate}

Define $\sub{d}(\sub{t})$ to be the element in ${\mathfrak{S}}_n$ that satisfies $\sub{t}=\sub{d}(\sub{t}) \sub{t}^{\bs{\lambda}}$.

\begin{defi} \rm   \label
Given a node $A=(1,c,d)$ the residue of $A$ is defined to be
\begin{equation} \label{residue bipartition}
\text{res}(A)= \left\{
        \begin{array}{ll}
          c-1+k, & \hbox{ if } d=1 \\
          c-1-k, & \hbox{ if } d=2
        \end{array}
      \right.
\end{equation}
where $k\in \Z$ such that $2k \equiv m \mod l $. For $\sub{t} \in \Std(n)$ and $1\leq j \leq n$, define the residue of $\sub{t} $ at $j$, as $r_{\sub{t}}(j):= \text{res}(A)$ where $A$ is the node occupied by $j$ in $\sub{t}$.
\end{defi}

Let $ \approx $ be the equivalence relation on $\Std(n)$ given by $\sub{s} \approx \sub{t}$ if
$r_{\sub{s}}(j)=r_{\sub{t}}(j) $ for $k =  1,2,\ldots ,n$.   The equivalence classes for $\approx$ are parameterized by elements in $I^{n}$, we call such elements residue sequences.
For $\bs{i} \in I^n$ we denote by $\Std(\bs{i})$ the corresponding class.
Any tableau $\sub{t}$ gives rise to a residue sequence that is denoted  by $\bs{i}^{\sub{t}}$. Then we have
$\sub{t} \in  \Std(\bs{i}^{\sub{t}})$, but in general $\Std(\bs{i})$  may be empty. In the particular case where $\sub{t} = \sub{t}^{\bs{\lambda}}$, for some $\bs{\lambda} \in \bi{n}$,  we write $\bs{i}^{\bs{\lambda}}$ instead of $\bs{i}^{\sub{t}^{\bs{\lambda}}}$.

\begin{defi}\label{bases}
Suppose that $\bs{\lambda} \in \bi{n}$ and $\sub{s}, \sub{t} \in \Std (\bs{\lambda})$.
Let $\de{s} =s_{i_1}\ldots s_{i_k}$ and $\de{t}=s_{j_1} \ldots s_{j_l}$ be reduced expressions for $\de{s}$ and $\de{t}$. Define
$$  \psi_{\sub{st}} :=    \psi_{i_1}\ldots \psi_{i_k}   e(\bs{i}^{\bs{\lambda}})
 \psi_{j_l}\ldots \psi_{j_1} \in \blob.     $$
\end{defi}

The elements $\psi_{\sub{st}}$ do not depend on the particular choices of a reduced expression for $\de{s}$ and $\de{t}$. Note that all elements defined above are homogeneous since they are products of homogeneous generators. Define the degree function $\deg: \Std(n)  \rightarrow \mathbb{Z}$ as
\begin{equation} \label{grado}
 \deg \sub{t} :=  \deg \psi_{\sub{tt}^{\bs{\lambda}}}.
\end{equation}
\begin{teo} \cite[Theorem 6.8]{Steen-David}\label{teo graded cellular bases}
The algebra $\blob$ is a graded cellular algebra with weight poset $(\bi{n}, \succeq)$,  $T(\bs{\lambda})= Std(\bs{\lambda})$ for $\bs{\lambda} \in \bi{n}$,  graded cellular basis
$\{ \psi_{\sub{st}}  \mid \sub{s,t} \in \Std(\bs{\lambda}) ; \bs{\lambda} \in \bi{n}  \}$ and degree function $\deg$ as in $(\ref{grado})$.
\end{teo}

Now that $\{\psi_{\sub{st}} \}$ is known to be a graded cellular basis we can define the graded cell and simple
$\blob$-modules which we denote by $\Delta(\bs{\lambda})$  and $L(\bs{\lambda})$, respectively. Thus the cell module has $\C$-basis
$\{ \psi_{\sub{t}} \mid \sub{t} \in \Std(\bs{\lambda}) \}$ and the $\blob$-action comes from the definition of cell modules. Although we do not need it, we remark that this cell modules coincide with the diagrammatic or standard modules for $\blob$. The main goal in this paper is to find the Laurent polynomials $[\Delta(\bs{\mu}):L(\bs{\lambda})]_t$ for all $\bs{\mu}, \bs{\lambda} \in \bi{n}$, we shall refer to these polynomials as graded decomposition numbers for $b_{n}(m)$.

\section{Degree function}
In this section we obtain a combinatorial formula for the degree function (\ref{grado}). We remark that in general this function does not coincide  with the degree function defined \cite[(3.5)]{bkw graded}, the main reason for this is that we do not work with the dominance order on $\Std(n)$ and
both degree functions depend heavily on the order considered on $\Std(n)$. First we give an interpretation of the degree function in terms of addable and removable nodes, similar to the one given in \cite[Definition 4.7]{hu-mathas} for standard tableaux using the dominance order. Finally, we give a formula for the degree function depending on \emph{walks and walls} on the Pascal triangle.

\medskip
Let $\bs{\lambda} \in \bi{n}$. The node $\alpha=(1,c,d)$ is called an addable node of $\bs{\lambda} $ if $\alpha \not \in \diagramof{\lambda}$ and $\diagramof{\lambda} \cup \{ \alpha \}$ is the diagram of a one-line bipartition of $n+1$. Similarly, $\rho \in \diagramof{\lambda} $ is called a removable node of $\bs{\lambda}$ if $\diagramof{\lambda} \backslash \{ \rho\}$ is the diagram of a one-line bipartition of $n-1$. Note that  any one-line bipartition has exactly two addable nodes. Furthermore, a one-line bipartition may have one or two removable nodes.

\medskip
Given two nodes $\alpha =(1,c_1,d_1)$ and $\beta =(1,c_2,d_2)$ then $\alpha $ is  said to be below  $\beta $ if $c_1>c_2$, or $c_1=c_2$, $d_1=1$ and $d_2=2$.
The concept of \emph{to be below} could have been defined in terms of $\sub{t}^{\la}$. In fact, given two nodes $\alpha$ and $\beta$ choose a  bipartition $\bs{\lambda}$ such that $\alpha , \beta \in [\bs{\lambda}]$. Then, the node $\alpha$ is below $\beta$ if and only if $\sub{t}^{\bs{\lambda}}(\alpha) > \sub{t}^{\bs{\lambda}}(\beta)$. Using the dominance order there is a similar interpretation of the concept of \emph{to be below} introduced in
\cite[Section 4]{hu-mathas}. Since $\sub{t}^{\la}$ does not coincide with the unique maximal bitableau for the dominance order, the two concepts does not coincide in general. 

\medskip
Let $\sub{t} \in \text{Std}(\bs{\lambda})$. For $k=1,\ldots, n$ let $\addable{k}{} $ be the set of all addable nodes of the bipartition $\text{ Shape}(\sub{t}_k)$ which are below of $\sub{t}^{-1}(k)$. Similarly, let $\removable{k}{} $ be the set of all removable  nodes of the bipartition $\text{ Shape}(\sub{t}_k)$ which are below of $\sub{t}^{-1}(k)$. Now define the sets $\addable{k}{m}$ and $\removable{k}{m}$ by

$$\begin{array}{c}
  \addable{k}{m}= \{  \alpha \in \addable{k}{} \mbox{ } | \mbox{ } \text{res}(\alpha) = r_{\sub{t}}(k) \} \\
             \\
  \removable{k}{m}= \{  \rho \in \removable{k}{} \mbox{ } | \mbox{ } \text{res}(\rho) = r_{\sub{t}}(k) \}
\end{array}
$$
It is easy to check that the sets $\addable{k}{m}$ and $\removable{k}{m}$ are empty or contain a single node, for all $\sub{t} \in \Std(\la)$ and  $1\leq k \leq n$. Let $g$ be the function defined by
$$ \begin{array}{rcl}
     g:\text{Std}(n) & \rightarrow & \mathbb{Z} \\
     \sub{t} & \rightarrow &  \sum_{k=1}^{n} \left(   |\addable{k}{m}|-|\removable{k}{m}|   \right)
   \end{array}$$

Using the above notation we can now give a characterization of the degree function. We need the following Lemma.

\begin{lem} \label{functiong}
Let $\bs{\lambda} \in \bi{n}$ and $\sub{s,t}\in \text{Std}(\bs{\lambda})$. Assume that $\sub{s} \succ \sub{t}$  and $s_{r}\sub{s}=\sub{t}$. Then
$$ g(\sub{t})-g(\sub{s})=\deg (\psi_r e(\bs{i}^{\sub{s}}))     $$
\end{lem}
\begin{dem}
First note that $\mathcal{A}_{\sub{s}}^{m}(k)=\mathcal{A}_{\sub{t}}^{m}(k)$ and $\mathcal{R}_{\sub{s}}^{m}(k)=\mathcal{R}_{\sub{t}}^{m}(k)$, for all $k \neq r,r+1$, since $s_{r}\sub{s}=\sub{t}$. Hence
 $$ g(\sub{t})-g(\sub{s})= \sum_{k=r}^{r+1} \left(   |\addable{k}{m}|-|\removable{k}{m}|   \right) - \sum_{k=r}^{r+1} \left(   |\mathcal{A}_{\sub{s}}^{m}(k)|-|\mathcal{R}_{\sub{s}}^{m}(k)|   \right)   $$

Set $e(\bs{i}^{\sub{s}})=(i_1,\ldots,i_r,i_{r+1},\ldots , i_n)$. We remark that $i_j= r_{\sub{s}}(j)$. The numbers involved in the above sums depends on  $i_r-i_{r+1}$ modulo $l$, so we  split the proof in four cases according to the followings options:

\begin{equation}
i_r-i_{r+1} \equiv \left\{
                     \begin{array}{r}
                       -1 \mod l \\
                       0 \mod l \\
                       1 \mod l   \\
                      \text{ otherwise. }
                     \end{array}
                   \right.
\end{equation}
We consider the first case, the remaining three cases follow in a similar way. Thus we assume that $i_r-i_{r+1}\equiv -1 \mod l$. Note that the node occupied by $r$ in $\sub{t}$ is below  the node occupied by $r+1$ in $\sub{t}$ since $\sub{s} \succ \sub{t}$. In this setting,  the values involved in the sums are

$$   \begin{array}{cccc}
       |\mathcal{A}_{\sub{s}}^{m}(r)|=0 & |\mathcal{A}_{\sub{s}}^{m}(r+1)|=0 &  |\mathcal{R}_{\sub{s}}^{m}(r)|=1 &   |\mathcal{R}_{\sub{s}}^{m}(r+1)| =0  \\
        & & & \\
       |\mathcal{A}_{\sub{t}}^{m}(r)|=0 & |\mathcal{A}_{\sub{t}}^{m}(r+1)| =0 &  |\mathcal{R}_{\sub{t}}^{m}(r)|=0  &   |\mathcal{R}_{\sub{t}}^{m}(r+1)|=0
     \end{array}
$$
Therefore $ g(\sub{t})-g(\sub{s})= 1= \deg (\psi_r e(\bs{i}^{\sub{s}}))  $, completing the proof in this case.
\end{dem}

\begin{exa}  
Assume that $l=5$ and $m=2$. Let $\la=((5),(1))\in \bi{6}$ and  $\sub{s}, \sub{t} \in \Std(\la)$ given by 
$$  \eje    $$
Note that $\sub{s} \succ \sub{t}$, $s_5\sub{s}=\sub{t}$ and $e(\bs{i}^{\sub{s}})=(1,2,3,4,4,0)$. We also have  
$\mathcal{R}_{\sub{s}}^{m}(5)=\{\sub{s}^{-1}(4) \}$,  $\mathcal{R}_{\sub{s}}^{m}(j)=\emptyset$  for $ j\neq 5$ and  $\mathcal{A}_{\sub{s}}^{m}(j)= \mathcal{R}_{\sub{t}}^{m}(j) =\mathcal{A}_{\sub{t}}^{m}(j)=\emptyset $ for all  $1\leq j\leq 6$.
Thus $g(\sub{s})=-1 $ and $g(\sub{t})=0$. Therefore, $g(\sub{t})-g(\sub{s})=1=\deg \psi_{5}e(\bs{i}^{\sub{s}})$.
\end{exa}

\begin{cor} \label{grado-functiong}
Let $\bs{\lambda} \in \bi{n}$ and $\sub{t} \in \text{Std}(\bs{\lambda})$. Then $g(\sub{t})=\deg(t)$.
\end{cor}
\begin{dem}
By  \cite[Lemma 4.15]{Steen-David} there is a sequence of one-line standard bitableaux
$$ \sub{t}= \sub{t}_0\prec \sub{t}_{1} \prec \ldots \prec \sub{t}_{k-1} \prec \sub{t}_k=\sub{t}^{\bs{\lambda}}  $$
such that $s_{i_{j}}\sub{t}_{j-1}=\sub{t}_j$ for $1\leq j \leq k$ and $\sub{d}(\sub{t})=s_{i_{1}} \ldots s_{i_{k}}$ is a reduced expression for $\sub{d}(\sub{t})$. Now using the above Lemma \ref{functiong} and the fact that $g(\sub{t}^{\bs{\lambda}})=0$ we have

\begin{equation*}
    g(\sub{t}) = \sum_{j=1}^{k} \left( g(\sub{t}_{j-1})-g(\sub{t}_j)  \right) = \sum_{j=1}^{k} \deg (\psi_{i_{j}} e(\bs{i}^{\sub{t}_{j}}))= \deg (\psi_{i_{1}} \ldots \psi_{i_{k}} e(\bs{i}^{\bs{\lambda}}))   = \deg (\sub{t})
\end{equation*}  \end{dem}

We now explain the bijection between $\text{Std}(n)$ and the set of all walks in the Pascal triangle with $n$ edges. Label vertices on the Pascal triangle by pairs of numbers giving level
(row) and weight (column):

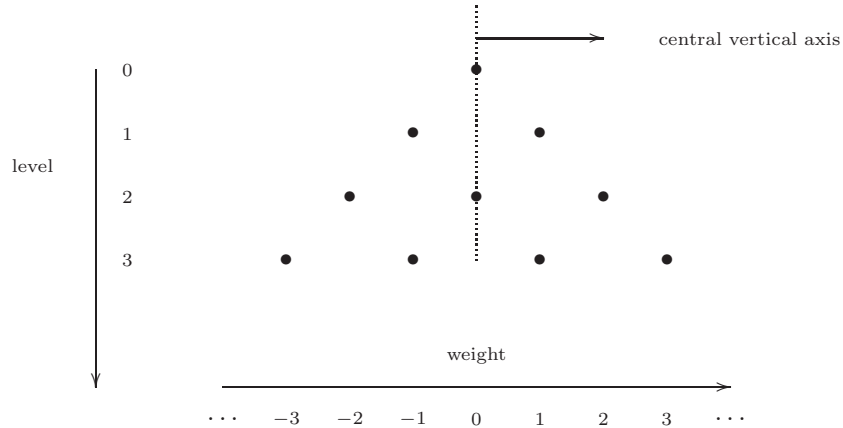
\[ \pascal  \]

\medskip
Label edges by vertex pairs: $((i, j), (i + 1, j \pm 1))$.  Let $W_n$ denote the set of walks from level 0 to level n on the Pascal triangle, that is graphs
with verices of the form $((i, j), (i+1, j \pm 1))$. There are various ways of specifying a particular $w\in W_n$.  First we can specify a walk as a sequence of
(level,weight) pairs. For example $((0, 0), (1, 1), (2, 0)))$. In particular, each edge is
a pair of such pairs. Since a walk on the Pascal triangle is always going down, we can even determine a walk simply as a sequence of weights, so our example becomes $(0,1,0)$. We use the above notation in the next section. Another way to describe walks is by a sequence of signs. We assign to each $((i, j), (i + 1, j - 1))$ edge the $-$ sign, while to each $((i, j), (i + 1, j + 1))$  edge we assign the $+$ sign. In this setting our example becomes $(+,-)$. With this at hand we can define a bijection between $\Std(n)$ and $W_n$.

\begin{defi} \rm  \label{tab-wal}
Let $\sub{t} \in \text{Std}(n)$. Define the bijective function $w$ by:
$$    \begin{array}{ccl}
        w: \text{Std}(n) & \rightarrow & W_n \\
        \sub{t} & \rightarrow & w(\sub{t})=(w(\sub{t})_1,\ldots, w(\sub{t})_n)
      \end{array}
$$
where $w(\sub{t})$ is the walk determined by the sign sequence

$$ w(\sub{t})_i= \left\{
                   \begin{array}{ll}
                     +, & \mbox{if } i \mbox{ is in the first component of } \sub{t}  \\
                     -, & \mbox{if } i \mbox{ is in the second component of } \sub{t}
                   \end{array}
                 \right.
   $$
\end{defi}

For $\sub{t} \in \text{Std}(n)$, when we specify the walk $w(\sub{t})$ as a sequence of weights we also write  $w(\sub{t})=(w(\sub{t})_0,w(\sub{t})_1,\ldots, w(\sub{t})_n)$, where in this setting the $w(\sub{t})_i$'s are integers. Note that for the sequence of signs there are $n$ signs, whereas for the sequence of weights there are $n+1$ integers, so there can be no confusion between the two notations.

\medskip
With this walk realization of the bitableaux, we can visualize the order $\succeq$. Indeed, let
 $\sub{s,t} \in \text{Std}(\bs{\lambda})$. Then $\sub{s} \succeq \sub{t}$ if and only if at each step of the two walks $w(\sub{s})$ is either strictly closer than $w(\sub{t})$
to the central vertical axis of the Pascal triangle or they are at the same distance
from the central axis and $w(\sub{s})$ is located (weakly) to the left of $w(\sub{t})$. In particular, for $\bs{\lambda}$, we can give a description of $w(\sub{t}^{\bs{\lambda}})$ as the walk that
first zigzags on and off the central vertical axis of the Pascal triangle, using the
signs $-$ and $+$ alternately, and then finishes using the sign $+$ repeatedly, if
$\bs{\lambda}$ is located in the positive half, or using the sign $-$ repeatedly, if
$\bs{\lambda}$ is located in the negative half. For example, the walk $w(\sub{t}^{\bs{\lambda}})$ for
$\bs{\lambda}= ((4),(2)) $ is:

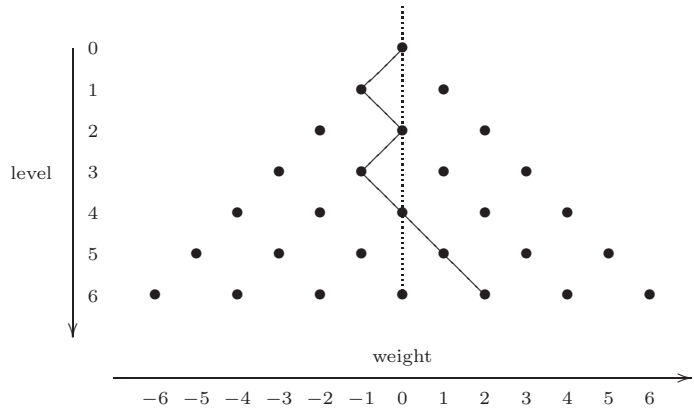
\[   \mayor \]

The walk realization of the bitableaux is also useful to express the residue sequence. In fact, for $\sub{s} \in \Std(n)$ the residue can be expressed in terms of the weight sequence of  $w(\sub{t})$ by the formula
\begin{equation}  \label{residuo-camino}
2r_{\sub{t}}(j)\equiv j-2+(w(\sub{t})_j-w(\sub{t})_{j-1})(w(\sub{t})_j+m ) \mod l
\end{equation}

\begin{lem} \label{grado-walks}
Let $\sub{t}\in \Std(n)$. Consider the walk $w(\sub{t})=(w(\sub{t})_0,w(\sub{t})_1, \ldots , w(\sub{t})_n)$ determined in Definition \ref{tab-wal} written as a sequence of weights. Define the sets
$$ A_{\sub{t}}^{1}=\{  1\leq j \leq n \mid   w(\sub{t})_{j}<0,  w(\sub{t})_{j-1} \equiv -m \mod l,  w(\sub{t})_{j} \equiv -m+1 \mod l    \} $$
$$ A_{\sub{t}}^{2}=\{  1\leq j \leq n \mid w(\sub{t})_{j}>0,   w(\sub{t})_{j-1} \equiv -m \mod l,    w(\sub{t})_{j} \equiv -m-1 \mod l   \}  $$
$$  R_{\sub{t}}^{1}=\{  1\leq j \leq n \mid w(\sub{t})_{j}<0,   w(\sub{t})_{j-1} \equiv -m-1 \mod l, w(\sub{t})_{j} \equiv -m \mod l   \}  $$
$$  R_{\sub{t}}^{2}=\{  1\leq j \leq n \mid w(\sub{t})_{j}>0,   w(\sub{t})_{j-1} \equiv -m+1 \mod l,  w(\sub{t})_{j} \equiv -m \mod l   \}  $$
Define also the sets $A_{\sub{t}}=A_{\sub{t}}^{1} \cup A_{\sub{t}}^{2}$ and $R_{\sub{t}}=R_{\sub{t}}^{1}\cup R_{\sub{t}}^{2}$. Then we have
\begin{description}
  \item[\rm (i)] $j\in A_{\sub{t}} $ if and only if $|\mathcal{A}_{\sub{t}}^{m}(j)|=1$
  \item[\rm (ii)] $j\in R_{\sub{t}} $ if and only if $|\mathcal{R}_{\sub{t}}^{m}(j)|=1$
\end{description}
\end{lem}
\begin{dem} We only prove (i), the result (ii) is proved similarly. Recall that $|\mathcal{A}_{\sub{t}}^{m}(j)|=0$ or $1$, for all $\sub{t} \in \Std(n)$ and $1\leq j \leq n$.  Suppose that $j \in A_{\sub{t}}$, then $j\in A_{\sub{t}}^{1}$ or $j\in A_{\sub{t}}^{2}$. Assume that $j\in A_{\sub{t}}^{1}$, the case $j\in A_{\sub{t}}^{2}$ is treated similarly. By definition of $A_{\sub{t}}^{1}$ we have
\begin{equation*} 
   w(\sub{t})_{j}<0, \qquad   w(\sub{t})_{j-1} \equiv -m \mod l, \qquad  w(\sub{t})_{j} \equiv -m+1 \mod l   
\end{equation*}
Recall that for any one-line bipartition there is two addable nodes, one in each component. Let $N_1$ and $N_2$ the addable nodes to $\Shape(\sub{t}_{j-1})$ in the first and second component, respectively. By $w(\sub{t})_{j-1} \equiv -m \mod l$ and (\ref{residuo-camino}), the nodes $N_1$ and $N_2$ have the same residue. On the other hand,
$ w(\sub{t})_{j} \equiv -m+1 \mod l $ implies that $N_1=\sub{t}^{-1}(j)$. Hence, the node $N_2$ is addable to  $\Shape(\sub{t}_{j})$ with the same residue of $N_1$. Finally, note that $N_2$ is below $N_1$ since $ w(\sub{t})_{j}<0$.  Therefore, $N_2 \in \mathcal{A}_{\sub{t}}^{m}(j)$. Consequently, $|\mathcal{A}_{\sub{t}}^{m}(j)|=1$.

\medskip
Conversely, suppose that $|\mathcal{A}_{\sub{t}}^{m}(j)|=1$. Let $N$ be the unique node in $\mathcal{A}_{\sub{t}}^{m}(j)$. Define $M$ to be the node occupied by $j$ in $\sub{t}$. Then, $N$ and $M$ have the same residue and $N$ is below $M$. 
Set $N= (1,c_1,d_1)$ and $M=(1,c_2,d_2)$. Using the conditions on the parameters $q$ and $m$ given in (\ref{qym}), one can check that two nodes with the same residue can not be located in the same column . Hence, $c_2\neq c_1$. Since $N$ is below $M$ we actually have  $c_2>c_1$. Note also that $N$ and $M$ can not be located in the same component. Hence, $d_1\neq d_2$. Assume that $d_2=1$, thus  $d_1=2$. Then, the fact that $N$ and $M$ have the same residue is equivalent to 
\begin{equation} \label{c1menosc2}
c_2-c_1 \equiv -m \mod l
\end{equation} 
On the other hand, recall that for any $1\leq i \leq n$ the weight $w(\sub{t})_i$ is equal to the number of nodes in the first component of $\Shape(\sub{t}_{|_i})$ minus the number of nodes in the second component of $\Shape(\sub{t}_{|_i})$. Therefore, it is easy to note that $w(\sub{t})_{j-1}= c_2-c_1$ and $w(\sub{t})_{j} = c_2-c_1+1 $. Then, by (\ref{c1menosc2}) and $c_2>c_1$ we obtain

\begin{equation*}
       w(\sub{t})_{j-1}  \equiv -m \mod l;     \qquad       w(\sub{t})_{j} \equiv -m+1 \mod l \qquad   \text{ and }    \qquad w(\sub{t})_{j} >0.
 \end{equation*}
This proves that $j\in A_{\sub{t}}^{1} \subset A_{\sub{t}}$. If $d_2=2$ then arguing similarly we get that $j\in A_{\sub{t}}^{2} \subset A_{\sub{t}}$, completing the proof.
\end{dem}

\begin{cor} \label{grado-sets A y R}
Let $\sub{t}\in \text{Std}(n)$. Let $A_{\sub{t}}$ and $R_{\sub{t}}$ the sets defined in the above Lemma. Then
\begin{equation}
\deg (\sub{t}) = |A_{\sub{t}}|-|R_{\sub{t}}|
\end{equation}
\end{cor}
\begin{dem} This is a direct consequence of Corollary \ref{grado-functiong} and Lemma \ref{grado-walks}.
\end{dem}

\medskip
We are now in position to give a graphical  interpretation on the Pascal triangle of the degree function, for this we need to draw \emph{walls} in the Pascal triangle. This means drawing vertical lines in each weight, $\bs{\lambda}$,  such that $\bs{\lambda}\equiv -m \mod l$, as shown in Figure for the case $l=5$ and $m=2$.

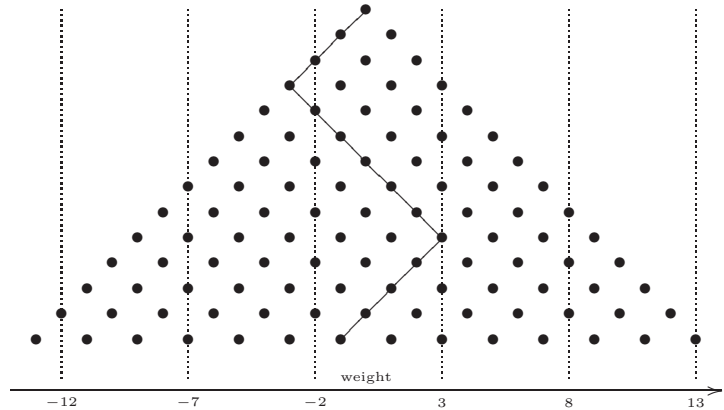
\[ \walls \]

Then, for $\sub{t} \in \Std(n)$, $|A_{\sub{t}}|$
(resp. $|R_{\sub{t}}|$) is the number of edges in the walk $w(\sub{t})$ such that the initial (resp. final) vertex is on a wall and the final vertex is closer than the initial vertex to the central axis of the Pascal triangle. For example, let $\sub{t}$ be the bitableau associated with the walk in the above figure, then  $A_t=\{  5,10 \} $ and $R_t=\{4\} $, consequently $\deg (\sub{t})=1$. It is also easy to check that $|A_{\sub{t}^{\bs{\lambda}}}|=|R_{\sub{t}^{\bs{\lambda}}}|=0 $  for all $\bs{\lambda} \in \bi{n}$, so $\deg (\sub{t}^{\bs{\lambda}})=0$.

\medskip
The walls drawn on the Pascal triangle define an alcove structure on $\mathbb{R}$, where the alcoves are the connected components of non-walls elements. We can thus refer to the alcove or wall in which a given weight lies. Note that by the conditions on the parameters (\ref{qym}), the weight $\bs{\lambda}=0$ always belongs to an alcove, that is, it is not on a wall. We refer to the alcove in which $\bs{\lambda} =0$ lies as the \emph{fundamental alcove}. Let $W$ be the infinite dihedral group on two generators $s_-$ and $s_+$, that is $W=\langle s_-,s_+ | s_{-}^{2}=s_{+}^{2}=1  \rangle$. The alcove structure defines an action of $W$ on $\mathbb{R}$,  by mapping $s_-$ (resp. $s_+$) to the reflection in the left  (resp. right) wall of the fundamental alcove.  Since the walls were drawn on integral weights, the subset $\Z$ of $\mathbb{R}$ is clearly invariant under this action.  Therefore, we can restrict  the action of $W$ to $\Z$. Let $\sim$ be the equivalence relation on $\Z$ determined by this action. The figure shows the orbit of $0$ under this action for $l=5$ and $m=2$. In this case, we have $\cdots \sim -10 \sim -4 \sim 0 \sim 6 \sim 10 \sim \cdots$.

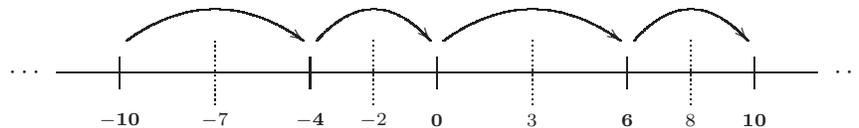
$$  \action $$

For a weight $\bs{\lambda} \in \Lambda_n$ we denote by $O_n(\bs{\lambda})$  the set of all $\bs{\mu} \in \Lambda_n$ such that $\bs{\mu} \sim \bs{\lambda}$.
Define $M_n(\bs{\lambda})$ to be the set
\begin{equation} \label{M- lambda}
M_n(\bs{\lambda})=\{ \bs{\mu} \in  O_{n}(\bs{\lambda}) \mid \text{ there exist } \sub{s} \approx \sub{t}^{\la} \text{ with } \Shape(\sub{s})=\bs{\mu} \}
\end{equation}

 Given a walk on the Pascal triangle we say that a subset of consecutive edges is
a \emph{wall to wall step} if these edges form a straight line between two walls of the same
alcove. A wall to wall step can be classified into three different types according to whether it crosses the fundamental alcove, it goes away from the central axis or it approaches the central axis. We denote by $F$, $O$ and $I$ to these types, respectively. For $\sub{s} \in \Std(n)$, we also define integers $\num{F}$, $\num{I}$ and $\num{O}$  as the number of occurrences  in $w(\mathfrak{s})$ of wall to wall steps of type $F$, $I$ and $O$, respectively. The following lemma is the first step in order to give an easy formula for the degree of $\sub{s} \in \Std(n)$ such that $\sub{s} \approx \sub{t}^{\bs{\lambda}}$, for some
$\bs{\lambda} \in \bi{n} $. Recall that  $\sub{s} \approx \sub{t}$ if  and only if $r_{\sub{s}}(j)=r_{\sub{t}}(j)$, for all $1\leq j \leq n$.

\begin{lem} \label{lema para grado positivo} Let $\bs{\lambda} \in \bi{n}$. A walk $w(\sub{s})$ for $\sub{s} \in \Std(n)$ satisfies $\sub{s} \approx \sub{t}^{\la}$ if and only if the following conditions hold:
\begin{description}
  \item[\rm (a)] First, $w(\sub{s})$ and $w(\sub{t}^{\bs{\lambda}})$ must match from level $0$ to the first contact of $w(\sub{t}^{\bs{\lambda}})$ with a wall of the fundamental alcove.
  \item[\rm (b)] Next, $w(\sub{s})$ makes wall to wall steps (as many as the number of alcoves between $\bs{\lambda}$ and $0$) of any type.
  \item[\rm (c)] Finally, $w(\sub{s})$ is completed with a straight line to the level $n$ in either direction.
\end{description}
\end{lem}
\begin{dem}
First, recall from (\ref{residuo-camino}) that for any $\sub{s} \in \Std(n)$ the residue can be expressed in terms of the weight sequence of  $w(\sub{s})$.
Assume that $\sub{s} \approx \sub{t}^{\bs{\lambda}}$ and suppose that we know the weights $w(\sub{s})_0,w(\sub{s})_1, \ldots , w(\sub{s})_{j-1}$ of the weight sequence of $w(\sub{s})$. Recall that $w(\sub{s})_{j}= w(\sub{s})_{j-1} \pm 1$, so $w(\sub{s})_{j}$ has only two options. Therefore, if $w(\sub{t})_{j-1} \not \equiv -m \mod l$ (equivalently, if $w(\sub{t})_{j-1}$ is not on a wall) only one option is acceptable to $w(\sub{s})_{j}$ because by replacing  in (\ref{residuo-camino}) we get two different values for $r_{\sub{s}}(j)$. But $r_{\sub{s}}(j)$ is a known value, actually $r_{\sub{s}}(j) =r_{\sub{t}^{\bs{\lambda}}}(j)$. Instead, if  $w(\sub{t})_{j-1}  \equiv -m \mod l$ (equivalently, if $w(\sub{t})_{j-1}$ is  on a wall) both options for  $w(\sub{t})_{j}$ are acceptable to $w(\sub{s})_{j}$ because by replacing  in (\ref{residuo-camino}) we get the same value for $r_{\sub{s}}(j)$. The lemma follows then by the description of $w(\sub{t}^{\bs{\lambda}})$ as the walk that first zigzags on and off the central vertical axis of the Pascal triangle, and then finishes with a straight line to the weight $\la$ at level $n$.
\end{dem}

\begin{rem} \rm
A walk as in the above lemma need not have parts (b) and (c). For example, if $\bs{\lambda}$ is located in the fundamental alcove, on one of its walls or in one of the two alcoves adjacent to the fundamental alcove a walk $w(\sub{s})$ with $\sub{s} \approx \sub{t}^{\bs{\lambda}}$ does not have part (b). On the other hand, if $\bs{\lambda}$ is located on a wall a walk $w(\sub{s})$ with $\sub{s} \approx \sub{t}^{\bs{\lambda}}$ does not have part (c). According to the above lemma, we split any walk $w(\sub{s})$ with $\sub{s} \approx \sub{t}^{\bs{\lambda}}$ in three parts (a), (b) and (c).
\end{rem}

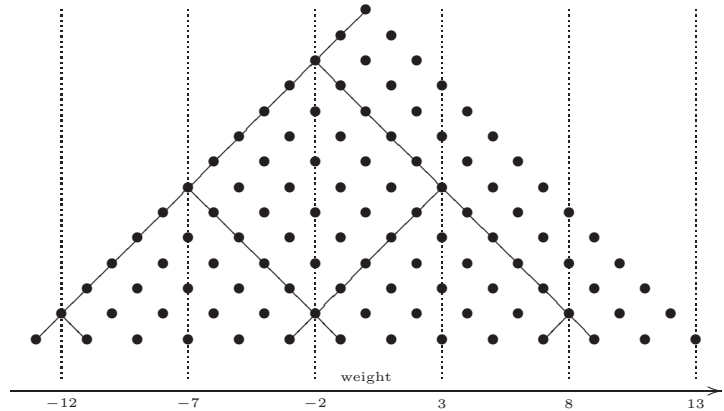
\[ \maximales \]

The above figure shows all walks on the Pascal triangle, $w(\sub{s})$, with $\sub{s} \approx \sub{t}^{\bs{\lambda}}$, where $\bs{\lambda}=((0),(13))\in \bi{13}$. In terms of the description given in the above lemma, for all walks in the figure, we have that  part (a) goes from level 0 to level 2,  part (b) goes from level $2$ to level $12$, and part (c) goes from level $12$ to level $13$. Hence, in this case $M_{13}(\la)=\{ -13, -11, -3,-1,7,9  \}$.

\begin{teo} \label{grado-positivo}
Let $\bs{\lambda} \in \bi{n}$ and $\sub{s} \in \text{Std}(n)$. Suppose that $\sub{s} \approx \sub{t}^{\bs{\lambda}}$. Then,
\begin{description}
  \item[\rm (a)] \label{grado positivo parte a} If $\sub{s} \in \text{Std}(\bs{\lambda})$ then $\sub{s}=\sub{t}^{\bs{\lambda}}$;
  %\item[\rm (b)]If $\bs{\mu}=\text{Shape}(\sub{s})  $ then $\bs{\mu} \in M_n(\bs{\lambda})$;
  \item[\rm (b)] If part (c) of $w(\sub{s})$ points toward the central axis  then $\deg (\sub{s})= \num{F} +1$;
  \item[\rm (c)] If part (c) of $w(\sub{s})$ does not points toward the central axis then $\deg (\sub{s})= \num{F}$.
\end{description}
\end{teo}
\begin{dem}
Statement (a) says that the unique $\sub{s} \in \Std(\la) $ with $\sub{s}\approx \sub{t}^{\la}$ is $\sub{t}^{\la}$, and this is clear from the Lemma \ref{lema para grado positivo}. Now, via Corollary \ref{grado-sets A y R}, we can conclude that the part (a) of  $w(\sub{s})$  has degree zero, in part (b) the degree is $\num{F}$, and  part (c) has degree $1$ or $0$ according to whether the final straight line points toward the central axis of the Pascal triangle or not. This proves (b) and (c).
\end{dem}

\section{Graded decomposition numbers}
In this section we obtain the main result in this paper, to determine the graded decomposition numbers for $b_n(m)$.  For $\bs{\lambda} \in \bi{n}$, denote by $b_{n}(m,\bs{\lambda})$ to be the subalgebra $e(\bs{i^{\lambda}})b_n(m) e(\bs{i^{\lambda}})$ of $b_n(m)$.   The basic strategy to find the graded decomposition numbers for $b_n(m)$ is to exploit Theorem \ref{gdn iguales teorema} on this subalgebra. We moreover need the known fact that the (ungraded) decomposition numbers of $b_{n}(m)$ are $0$ or $1$ (See for instance \cite[Theorem 5.5]{Steen1}). 

\begin{rem}   \rm
If we assume that the graded decomposition numbers for $b_n$ are polynomials with constant coefficient equal to zero then we can obtain these numbers using analogous methods to those used by Kleshchev and  Nash in \cite{nash}, without using the prior knowledge about the (ungraded) decomposition numbers for $b_n$ mentioned in the last paragraph. The additional hypothesis on the graded decomposition numbers can be proved by brute force calculations over the homogeneous presentation for $b_n$. However, for the sake of readability, we prefer the presentation of the article as it stands. 
\end{rem}

\begin{teo} \label{gcb subalgebra}For $\bs{\lambda} \in \bi{n}$,  we have that $b_{n}(m,\bs{\lambda})$ is a positively graded cellular algebra with
weight poset $(M_n(\bs{\lambda}), \succeq)$, $T(\bs{\mu})= \text{Std}(\bs{\mu}) \cap \text{Std}(\bs{i^{\lambda}}) $  for $\bs{\mu} \in M_n(\bs{\lambda})$, graded cellular basis
\begin{equation}
\{   \psi_{\sub{st}} \mbox{ } | \mbox{ } \sub{s,t} \in \text{Std}(\bs{\mu}) \cap \text{Std}(\bs{i^{\lambda}}) \text{ for }  \bs{\mu} \in  M_n(\bs{\lambda})  \}
\end{equation}
and degree function as (\ref{grado}).
\end{teo}
\begin{dem} By Theorem \ref{teo graded cellular bases} and the orthogonality of the KLR-idempotents, $e(\bs{i})$, we have that $b_{n}(m,\bs{\lambda})$  has a $\C$-basis consisting of all elements $\psi_{\sub{st}}$ such that
\begin{equation} \label{base subalgebra}
 \bs{i}^{\sub{s}}=\bs{i}^{\sub{t}}= e(\bs{i^{\lambda}})
\end{equation}
Thus, $b_{n}(m,\bs{\lambda})$  is a positively $\Z$-graded algebra since all bitableaux satisfying the condition (\ref{base subalgebra}) have non-negative degree by Theorem \ref{grado-positivo}. All claims about the cellularity of $b_{n}(m,\bs{\lambda})$ follow from the cellularity of $b_{n}(m)$ and the definition of $M_{n}(\la)$.\end{dem}

\medskip
Now that $\{   \psi_{\sub{st}} \mbox{ } | \mbox{ } \sub{s,t} \in \text{Std}(\bs{\mu}) \cap \text{Std}(\bs{i^{\lambda}}) \text{ for }  \bs{\mu} \in  M_n(\bs{\lambda})  \}$ is known to be a graded cellular basis for $b_{n}(m,\bs{\lambda})$ we can define graded cell and simple modules which we denote by $\Delta_{\bs{\lambda}}(\bs{\mu})$ and $L_{\bs{\lambda}}(\bs{\mu})$ respectively, for $\bs{\mu} \in  M(\bs{\lambda}) $. Note that
$e(\bs{i^{\lambda}})\Delta(\bs{\mu})=\Delta_{\bs{\lambda}}(\bs{\mu}) $ and $e(\bs{i^{\lambda}})L(\bs{\mu})=L_{\bs{\lambda}}(\bs{\mu})  $. We can also define graded decomposition numbers for $b_n(m,\la)$. Graded decomposition numbers of $\blob$ and $b_n(m,\la)$ are related via equation (\ref{gdn iguales}). 

\begin{lem}  \label{gdn zero}Let $\bs{\lambda} \in \bi{n}$. If $\bs{\mu} \not \in M_n(\bs{\lambda})$ then $[\Delta(\bs{\mu}):L(\bs{\lambda})]_t=0$.
\end{lem}
\begin{dem}  First, recall that for any $\bs{\nu} \in \bi{n}$ the cell $b_{n}$-module $\Delta(\bs{\nu})$ has a $\C$-basis
$\{ \psi_{\sub{s}} \mbox{  } | \mbox{  } \sub{s} \in \text{Std} (\bs{\nu})  \}$ and that by the orthogonality of the KLR-idempotents, $e(\bs{i})$, these act on this basis according to the rule
\begin{equation}
e(\bs{i})\psi_{\sub{s}} = \left\{
                            \begin{array}{rl}
                              \psi_{\sub{s}}, & \text{ if } \bs{i}^{\sub{s}}=\bs{i} \\
                                            0,& \text{ if } \bs{i}^{\sub{s}}\neq \bs{i}
                            \end{array}
                          \right.
\end{equation}
Now, if $\bs{\mu} \not \in M_n(\bs{\lambda})$ then for all $\sub{s} \in \Std(\bs{\mu})$ we have $\bs{i}^{\sub{s}}\neq \bs{i}^{\la}$, thus
\begin{equation*}
e(\bs{i}^{\la})\Delta (\bs{\mu})=\Delta_{\bs{\lambda}}(\bs{\mu}) = 0
\end{equation*}
By the description given in Lemma \ref{lema para grado positivo} for all walks $w(\sub{s})$ with\ $\bs{i}^{\sub{s}}=\bs{i}^{\la}$ (and therefore for all $\bs{\mu} \in M_{n}(\bs{\la})$), it is straightforward to check that  $\la$ is a minimal element of $M_n(\la)$. Furthermore, Theorem \ref{grado-positivo}(a) implies that $\sub{t}^{\bs{\lambda}}$ is the unique standard bitableau in $ \Std(\la)$ with residue sequence equal to $\bs{i}^{\la}$.  Hence
\begin{equation*}
e(\bs{i}^{\la})\Delta (\la)=\Delta_{\la}(\la) = L_{\la}(\la) = \text{Span}_{\C}\{  \psi_{\sub{t}^{\bs{\lambda}}} \}
\end{equation*}
Therefore, by (\ref{gdn iguales}) we have
\begin{equation*}
 [\Delta(\bs{\mu}):L(\bs{\lambda})]_t= [\Delta_{\la}(\bs{\mu}):L_{\la}(\la)]_t=0
\end{equation*}
completing the proof of the lemma.
\end{dem}

\subsection{The non-wall case.}
If $\la \in \bi{n}$ belongs to the fundamental alcove it is straightforward to check that $M_{n}(\bs{\lambda})= \{ \bs{\lambda} \}$. Hence, by Lemma \ref{gdn zero} the module $L(\bs{\lambda})$ only appears as a graded composition factor in $\Delta(\bs{\lambda})$, and in this case by Theorem \ref{gdn in a cellular algebra}(c)
 $ [\Delta(\bs{\lambda}):L(\bs{\lambda}) ]_t=1 $. Therefore, we fix $\bs{\lambda} \in \bi{n}$ that does not belong to the fundamental alcove. Furthermore, for the rest of this subsection we also assume that $\bs{\lambda}$ is not on a wall of the Pascal triangle. The other case  will be treated in the forthcoming subsection.

\medskip
Since $b_n(m,\la)$ is a positively graded cellular algebra  $\Delta_{\la}(\bs{\mu})$ is also positively graded, for all $\bs{\mu} \in M_n(\la)$. Then 
$\dim_{t} \Delta_{\la}(\bs{\mu})\in \Z[t]$, for all $\bs{\mu} \in M_n(\la)$. Again by the positive grading on $b_n(m,\la)$ we can define the grading filtration for each $\Delta_{\la}(\bs{\mu})$. In order to know the dimensions of the quotients that appear in this grading filtration it is enough with to know the coefficients of $\dim_{t} \Delta_{\la}(\bs{\mu}) \Z[t]$. This is our next goal. We derive the graded decomposition number for $b_n$ from this. For $\bs{\lambda} \in \bi{n}$ define the number $\kappa(\bs{\lambda}) $ as the number of alcoves between $\bs{\lambda}$ and $0$.

\begin{lem} \label{cantidad M_n}
Let $\bs{\lambda} \in \bi{n}$. Then $|M_n(\bs{\lambda})|= 2(\kappa(\bs{\lambda})+1)$.
\end{lem}
\begin{dem}
Since $\bs{\lambda}$ is not on a wall, there is in each alcove a unique representative for the orbit of $\bs{\lambda}$. So, if $\bs{\lambda}$ is located in the positive (resp. negative) side of the Pascal triangle then to the right (resp. left)  of the (resp. left) right wall of the fundamental alcove there is exactly $\kappa(\bs{\lambda})+1$ elements in $M_n(\bs{\lambda})$. Reflecting these elements through the right (resp. left) wall of the fundamental alcove we can get all the elements in $M_n(\bs{\lambda})$ to the left (resp. right) of such wall. Hence, $|M_n(\bs{\lambda})|= 2(\kappa(\bs{\lambda})+1)$.
\end{dem}

\medskip
In order to give a precise description of $ \dim_t(\Delta_{\bs{\lambda}}(\bs{\mu}) )$, for $\bs{\mu} \in M_n(\bs{\lambda})$, we need to index the set  $M_n(\bs{\lambda})$. Set $\bs{\lambda}_1=\bs{\lambda}$. Assume that $\bs{\lambda}$ is located in the negative (resp. positive) side of the Pascal triangle. For
$1\leq i <|M_n(\bs{\lambda})|$ and $i$ odd, define $\bs{\lambda}_{i+1}$ as the rightmost  (resp. leftmost) weight in $M_{n}(\bs{\lambda}) \backslash \{ \bs{\lambda}_{j}\}_{j=1}^{i}$. On the other hand, if $1\leq i <|M_n(\bs{\lambda})|$ and $i$ even then we define  $\bs{\lambda}_{i+1}$ as the leftmost (resp. rightmost) weight in $M_{n}(\bs{\lambda}) \backslash \{ \bs{\lambda}_{j}\}_{j=1}^{i}$.

\begin{exa} \rm \label{ejemplo}
In the figure is shown the indexation of $M_n(\bs{\lambda})$, for $\bs{\lambda}=-19$, $l=5$ and $m=2$.
$$ \indexation  $$
\end{exa}

\begin{lem} \label{ayuda para el importante} Let $\bs{\lambda} \in \bi{n}$ and $\bs{\lambda}_{4j+1} \in M_{n}(\bs{\lambda}) $. For $0\leq i \leq j $, define
\begin{equation} \label{conjunto D_subi}
D_{i}^{j}= \{  \sub{s} \in \Std(\bs{\lambda}_{4j+1}) \mbox{  } | \mbox{  }  \sub{s} \approx \sub{t}^{\bs{\lambda}} \text{ and } \deg(s)=  2i   \}
 \end{equation}
Then, all $\sub{s} \in \Std(\bs{\lambda}_{4j+1})$ with $\sub{s} \approx \sub{t}^{\bs{\lambda}}$  belong to some  $D_i^{j}$ and
$|D_{i}^{j}|=| \Std(\mu_{i})|$, where $\mu_{i}^{j}$ is the two-column partition of $\kappa (\bs{\lambda})$ given by
\begin{equation}  \label{partition associated}
\mu_{i}^{j}=(\kappa(\bs{\lambda}) -j+i, j-i)'
\end{equation}
\end{lem}

\begin{dem}
Let $\bs{\lambda}_{4j+1} \in M_{n}(\bs{\lambda})$. Let $\sub{s} \in \Std(\bs{\lambda}_{4j+1})$ with $\sub{s} \approx \sub{t}^{\bs{\lambda}}$. According to Lemma \ref{lema para grado positivo} we can split the walk $w(\sub{s})$ in three parts (a), (b) and (c). Now, by a routine analysis of the indexation given for  $M_n(\bs{\lambda})$,  it is clear that $\bs{\lambda}_{4j+1}$  and $\bs{\lambda}$ are on the same side (positive or negative) of the Pascal triangle (actually, all weights in $M_{n}(\bs{\lambda})$ with odd subscript are on the same side),  $\bs{\lambda}_{4j+1}$ does not belong to the fundamental alcove, and $\bs{\lambda}_{4j+1}$  is located  $2j$ alcoves closer than $\bs{\lambda}$ to the fundamental alcove.

\medskip
Furthermore, part (c) of $w(\sub{s})$  always points away from the central axis of Pascal triangle (actually, this line coincide for all $\sub{s}$ under the above conditions). Thus,  Theorem \ref{grado-positivo}(c) implies that $\deg (\sub{s}) = n_{\sub{s}}(F)$, where we recall that $n_{\sub{s}}(F)$ was defined as the number of occurrences  in $w(\mathfrak{s})$ of wall to wall steps of type $F$ (similarly, we have defined the integers $n_{\sub{s}}(I)$ and $n_{\sub{s}}(O)$). Note also that  $n_{\sub{s}}(F)$ is even, because $\bs{\lambda}_{4j+1}$  and $\bs{\lambda}$ are on the same side  of the Pascal triangle, and for  $\sub{s} \in \Std(\bs{\lambda}_{4j+1}) $ we have  $0\leq n_{\sub{s}}(F) \leq 2j $ because $\bs{\lambda}_{4j+1}$  is located  $2j$ alcoves closer than $\bs{\lambda}$ of the fundamental alcove. Consequently, $\deg (\sub{s})$ is even and $0\leq \deg(\sub{s}) \leq 2j $. This proves the first claim of the Lemma.

\medskip
Let $0\leq i \leq j$ and  assume that $\sub{s} \in D_{i}^{j}$. By the above paragraph,  we can replace the condition $\deg (\sub{s})= 2i$ in the definition of $D_{i}^{j}$ by $n_{\sub{s}}(F)=2i$. Parts (a) and (c) of $w(\sub{s})$ are fixed, and part (b) (and therefore, the entire walk $w(\sub{s})$) is determined by a sequence of wall to wall steps. More conveniently, we can describe the walk $w(\sub{s})$ as an ordered word in the alphabet with three letters $\{ F,I,O \}$ in the obvious way. For example, for
 $\bs{\lambda} =((2),(20)) \in \bi{22}$, $l=5$, and $m=2$ we have that the bitableau described in the figure below as a walk on the Pascal triangle is in $\text{Std}(\bs{\lambda}_5)$, where $\bs{\lambda}_5 =((7),(15))$ and in terms of ordered words on $\{ F,I,O \}$ correspond to $FFO$.

$$ \dibujo $$

We now associate to  $\sub{s} $ a two-column standard tableau of shape $\mu_{i}^{j}$ as follows: for $w(\mathfrak{s})$  described as a ordered word in $\{F,O,I\}$, the tableau associated to $\sub{s}$ is  determined by placing on the second column the entries corresponding to the positions at which the letter $I$ appears in the respective ordered word. Therefore, the shape of the two-column partition associated to $\sub{s}$ is $(\num{F}+\num{O},\num{I})'$.

\medskip
In order to check that the above assignment is well defined first note that at the first $k$ positions of the ordered word associated to $w(\sub{s})$ the number of occurrences of the letter $O$ is greater or equal than the number of occurrences of the letter $I$, this shows that the two-column tableau assigned to $\sub{s}$ is standard. Next, recall that $\kappa(\bs{\lambda})$ is the number the alcoves between $\bs{\lambda}$ and $0$, so by
Lemma \ref{lema para grado positivo}(b) we have
\begin{equation} \label{ayuda1}
\kappa(\bs{\lambda})=n_{\sub{s}}(F)+n_{\sub{s}}(O)+n_{\sub{s}}(I) = 2i +n_{\sub{s}}(O)+n_{\sub{s}}(I)
\end{equation}
On the other hand, since $\bs{\lambda}_{4j+1}$ is located  $2j$ alcoves closer than $\bs{\lambda}$ to the fundamental alcove we have
\begin{equation}\label{ayuda2}
\kappa(\bs{\lambda})= 2j +n_{\sub{s}}(O)-n_{\sub{s}}(I)
\end{equation}
Combining (\ref{ayuda1}) and (\ref{ayuda2}) we obtain $n_{\sub{s}}(I)=j-i$. This proves that the two-column tableau associated to $w(\sub{s})$ has actually shape $\mu_{i}^{j}$. Therefore, the assignment is well defined. Finally, for any two-column standard tableau of shape $\mu_{i}^{j}$ is straightforward to check that one can recover a walk $w(\sub{s})$ with $\sub{s} \in D_{i}^{j}$. Hence, the above assignment is a bijection and $|D_{i}^{j}| =|\Std(\mu_{i}^{j})|$, completing the proof of the Lemma.
\end{dem}

\begin{teo}  \label{gd standard}
Let $\bs{\lambda} \in \bi{n}$ and assume that $\kappa(\bs{\lambda})\geq 1$. The graded dimension of $\Delta_{\bs{\lambda}}(\bs{\lambda}_i)$, for
$ \bs{\lambda}_i \in M_n(\bs{\lambda}) $, is completely determined by the formulas:

$$ \begin{array}{l}
 (a)  \qquad   \dim_t(\Delta_{\bs{\lambda}}(\bs{\lambda}_{4j+1}) )= \sum_{i=0}^{j} c_{i}t^{2i} \\
  %\qquad \text{ for } 0\leq 4j \leq \left[ \frac{2\kappa(\bs{\lambda}) +1}{4}   \right]
   \\
(b) \qquad \dim_t(\Delta_{\bs{\lambda}}(\bs{\lambda}_{4j+2}) )= t \dim_t(\Delta_{\bs{\lambda}}(\bs{\lambda}_{4j+1}) )     \\
      \\
(c) \qquad \dim_t(\Delta_{\bs{\lambda}}(\bs{\lambda}_{4j+3}) )= t \dim_t(\Delta_{\bs{\lambda}}(\bs{\lambda}_{4j+1}) )     \\
      \\
 (d) \qquad  \dim_t(\Delta_{\bs{\lambda}}(\bs{\lambda}_{4j+4}) )=t^{2} \dim_t(\Delta_{\bs{\lambda}}(\bs{\lambda}_{4j+1}) )
   \end{array}
$$
where $c_i=|\Std(\mu_{i}^{j})|$ and $\mu_{i}^{j}$ is the two-column partition of $\kappa(\bs{\lambda})$ defined in (\ref{partition associated}).
\end{teo}

\begin{dem} By Theorem \ref{gcb subalgebra} we have
\begin{equation}  \label{ecuacion de dimension graduada}
\dim_t(\Delta_{\bs{\lambda}}(\bs{\lambda}_i) )= \sum_{\substack{ \mathfrak{s}\in \Std(\bs{\lambda}_i) \\ \bs{i}^{\mathfrak{s}}=\bs{i^{\lambda}}}} t^{\deg (\mathfrak{s})}
\end{equation}
Hence, part (a) in the Theorem follows immediately from Lemma \ref{ayuda para el importante}. Now, we prove (b). Assume that $\bs{\lambda}$ is located on the negative (resp. positive) side of the Pascal triangle. Then, by the indexing on $M_{n}(\bs{\lambda})$, the weight $\bs{\lambda}_{4j+2}$ is obtained from $\bs{\lambda}_{4j+1}$ by reflection about the left (resp. right) wall of the fundamental alcove (See Example \ref{ejemplo}). Next, define the sets
$$  A=\{  \sub{s} \in \Std(\bs{\lambda}_{4j+1}) \mbox{ } |\mbox{ } \sub{s} \approx \sub{t}^{\bs{\lambda}}  \}    \qquad  B=\{  \tilde{\sub{s}} \in \Std(\bs{\lambda}_{4j+2})\mbox{ } | \mbox{ } \tilde{\sub{s}} \approx \sub{t}^{\bs{\lambda}}  \}    $$
For $\sub{s} \in A$ we associate an element $\tilde{\sub{s}}\in B$ in the following way: Let $l$ be the level at which the walk $w(\sub{s})$
intersects the left (resp. right) wall of the fundamental alcove for the last time. Then the walks $w(\tilde{\sub{s}})$  and $w(\sub{s})$  match from level $0$ to level $l$, and then from level $l$ to level $n$, $w(\tilde{\sub{s}})$ is obtained from $w(\sub{s})$ by reflection about the left (resp. right) wall of the fundamental alcove. It is clear that this process is reversible, so it defines a bijection between $A$ and $B$, and via Theorem \ref{grado-positivo}
we can conclude that  $\deg(\tilde{\sub{s}})= \deg(\sub{s})+1$. Consequently, by (\ref{ecuacion de dimension graduada}) we have
$$
t\dim_t(\Delta_{\bs{\lambda}}(\bs{\lambda}_{4j+1}) ) = t\sum_{\sub{s} \in A} t^{\deg (\mathfrak{s})}  \\
                                                    = \sum_{\sub{s} \in A} t^{\deg (\mathfrak{s})+1}  \\
                                                    = \sum_{\tilde{\sub{s}} \in B} t^{\deg (\tilde{\mathfrak{s}})}  \\
                                                     = \dim_t(\Delta_{\bs{\lambda}}(\bs{\lambda}_{4j+2}) )
$$
proving (b). Parts (c) and (d), follow in a similar way.
\end{dem}

\begin{cor} Let $\bs{\lambda} \in \bi{n}$ and suppose that $\kappa(\bs{\lambda})\geq 1$.  For $r=1,2,3,4$ and $\bs{\lambda}_{4j+r} \in M_{n}(\bs{\lambda})$ we have
\begin{equation} \label{dimension standard}
\dim_{\C}(\Delta_{\bs{\lambda}} (\bs{\lambda}_{4j+r}))= \sum_{i=0}^{j}|Std(\mu_{i}^{j})|
\end{equation}
where $\mu_{i}^{j}$ is the two-column partition of $k(\bs{\lambda})$ defined in (\ref{partition associated}).
\end{cor}
\begin{dem} This follows immediately by putting  $t=1$ in the above Theorem.
\end{dem}

\begin{cor} \label{casi simple modules} Let $\bs{\lambda} \in \bi{n}$ with $\kappa(\bs{\lambda})\geq 1$. Then, $L_{\bs{\lambda}}(\bs{\lambda}_k)\neq 0$ if and only if $k=4j+1$. Furthermore,
\begin{equation}  \label{dimension de los simples}
\dim_{\C} L_{\bs{\lambda}}(\bs{\lambda}_{4j+1})=\dim_{t} L_{\bs{\lambda}}(\bs{\lambda}_{4j+1}) = |\Std(\mu_{0}^{j})|
\end{equation}
where $\dim_{\C} L_{\la}(\la_{4j+1})$ is viewed as polynomial over $t$ in the natural way.
\end{cor}
\begin{dem} Since $b_n(m,\bs{\lambda})$ is a positively graded cellular algebra the modules $L_{\bs{\lambda}}(\bs{\lambda}_k)$ are pure of degree zero.  Therefore,   $\dim_{\C} L_{\bs{\lambda}}(\bs{\lambda}_{4j+1})=\dim_{t} L_{\bs{\lambda}}(\bs{\lambda}_{4j+1})$ and using  Theorem \ref{gd standard}  we can also conclude that
 \begin{equation} \label{acotacion dimension simples}
\dim_{\C} (L_{\bs{\lambda}}(\bs{\lambda}_{k})) \leq \dim_{t} \Delta_{\bs{\lambda}}(\bs{\lambda}_{k})_{t=0} = \left\{
                                                                                                                           \begin{array}{ll}
                                                                                                  |\Std(\mu_{0}^{j})| & \text{if } k=4j+1 \\
                                                                                                                  0         & \text{otherwise}.
                                                                                                                           \end{array}
                                                                                                                         \right.
\end{equation}
Then, $L_{\bs{\lambda}}(\bs{\lambda}_k)\neq 0$  only if $k=4j+1$.  On the other hand, recall that for $b_n(m)$ is known that the (ungraded) decomposition numbers are $0$ or $1$, therefore by putting $t=1$ in (\ref{gdn iguales}) this is also true for $b_n(m,\bs{\lambda})$. Now note that $\bs{\lambda}_{|M_{n}(\bs{\lambda})|}$ is in the fundamental alcove, so this is the maximal element in $M_{n}(\bs{\lambda})$ with respect to the order $\succeq$. Hence, by (\ref{dimension standard}) and (\ref{acotacion dimension simples})
\begin{align*}
\text{dim}_{\C} \Delta_{\bs{\lambda}}(\bs{\lambda}_{|M_{n}(\bs{\lambda})|}) & \leq \sum_{\scriptscriptstyle{ \bs{\lambda}_{4j+1} \in M_{n}(\bs{\lambda})}} \text{dim}_{\C} L_{\bs{\lambda}}(\bs{\lambda}_{4j+1}) \\
                                                                            & \leq\sum_{\scriptscriptstyle{ \bs{\lambda}_{4j+1} \in M_{n}(\bs{\lambda})}}   |\Std(\mu_0^{j})| \\
                                                                            & \leq \sum_{ \scriptscriptstyle{ \mu \in \pacol{\kappa(\bs{\lambda})} } }    |\Std(\mu)|  \\
                                                                           &= \text{dim}_{\C} \Delta_{\bs{\lambda}}(\bs{\lambda}_{|M_{n}(\bs{\lambda})|})
\end{align*}
Therefore, the inequalities  become equalities and $\text{dim}_{\C} (L_{\bs{\lambda}}(\bs{\lambda}_{4j+1})) = |\Std(\mu_0^{j})|$. Thus  $L_{\bs{\lambda}}(\bs{\lambda}_k)\neq 0$ if $k=4j+1$.
\end{dem}

%   9082948-3

\begin{rem} \label{dimension simples mayor que uno} \rm
Assume that $\kappa(\bs{\lambda}) \geq 1$ and $\kappa(\bs{\lambda})\neq 2$. Then, we have
\begin{equation*}
|\Std(\mu_0^{j})|  \left\{
                            \begin{array}{ll}
                              =1, & \text{if } j=0 \\
                              >1, & \text{if } j\neq0
                            \end{array}
                          \right.
\end{equation*}
Thus, under the above conditions on  $\kappa(\bs{\lambda})$, the algebra $b_n(m,\bs{\lambda})$  has a unique (up degree shift) one-dimensional graded simple module. Such module is $L_{\bs{\lambda}}(\bs{\lambda}_{4j+1})$ when $j=0$, that is, $L_{\bs{\lambda}}(\bs{\lambda}_{1})=L_{\bs{\lambda}}(\bs{\lambda})$ since $\bs{\lambda}_1=\bs{\lambda}$. If $\kappa(\bs{\lambda})=2$ then  by  Lemma \ref{cantidad M_n} we have $|M_{n}(\bs{\lambda})|=6$. Hence, by the above Corollary $b_n(m,\bs{\lambda})$ has two (up degree shift) non-isomorphic simple modules $L_{\bs{\lambda}}(\bs{\lambda}_{1})$ and $L_{\bs{\lambda}}(\bs{\lambda}_{5})$, both of dimension one.
\end{rem}

\medskip
We are now able to prove the main theorem in this paper for the \emph{non-wall case}, to provide the graded decomposition numbers for the blob algebra.  Surprisingly, it was more difficult to determine the graded decomposition numbers for the case $\kappa(\bs{\lambda})=2$ than for the general case. This difficulty lies in the fact that for $\kappa(\bs{\lambda})= 2$   there is two one-dimensional simple module for $b_n(m,\bs{\lambda})$, as explained in the previous remark.

\begin{teo} \label{main theorem} Let $\bs{\lambda} \in \bi{n}$. For $\bs{\lambda}_{k}\in M_{n}(\bs{\lambda})$  we have
 \begin{equation} \label{gdn answer}
 [\Delta(\bs{\lambda}_k) : L(\bs{\lambda})]_t= \left\{
                                          \begin{array}{ll}
                                            t^{2j},   & \text{ if } k=4j+1;  \\
                                            t^{2j+1}, & \text{ if } k=4j+2;  \\
                                            t^{2j+1}, & \text{ if } k=4j+3;  \\
                                            t^{2j+2}, & \text{ if } k=4j+4.
                                          \end{array}
                                        \right.
 \end{equation}
\end{teo}
\begin{dem}
By Remark \ref{dimension simples mayor que uno} we know that $\dim_{t}L_{\la}(\la)=1$ so Theorem \ref{gdn iguales} implies
\begin{equation*}
 [\Delta(\bs{\lambda}_k) : L(\bs{\lambda})]_t=  [\Delta_{\la}(\bs{\lambda}_k) : L_{\la}(\bs{\lambda})]_t
\end{equation*}
Therefore, we prove the theorem for the graded decomposition numbers of $b_n(m,\la)$, $[\Delta_{\la}(\bs{\lambda}_k) : L_{\la}(\bs{\lambda})]_t$. On the other hand, Theorem \ref{gdn in a cellular algebra}(d)  relates the graded dimension of cell and simple modules with the graded decomposition numbers via the formula
\begin{equation}  \label{relacion dimension cell y simples}
\dim_{t}\Delta_{\la}(\la_{k})= \sum_{\la_{j} \preceq \la_{k}} [\Delta_{\la}(\la_{k}):L_{\la}(\la_{j})]_t  \dim_{t} L_{\la} (\la_{j})
\end{equation}
Assume that $\kappa(\bs{\lambda})=0$, then $\la$ is located in one of the two alcoves adjacent to the fundamental alcove, and $|M_n(\la)|=2$ by Theorem \ref{cantidad M_n}. Write $M_{n}(\la)=\{\la_1,\la_2\}$. Then, $\la_1=\la$ and $\la_2$ is in the fundamental alcove. Combining Lemma \ref{lema para grado positivo} and Theorem \ref{grado-positivo} we have
\begin{equation} \label{kapalambda 0}
  \dim_t\Delta_{\la}(\la_2)=t
\end{equation}
Since $b_n(m,\la)$ is a positively graded cellular algebra the modules $L_{\la}(\la_k)$, $k=1,2$  are pure of degree zero. Then, $L_{\la}(\la_2)= 0$ by (\ref{kapalambda 0}). By Theorem \ref{gdn in a cellular algebra} we know that $[\Delta_{\la}(\la): L_{\la}(\la)]_t=1 $ and for $k=2$ equation (\ref{relacion dimension cell y simples}) becomes $t=[\Delta_{\la}(\la_2): L_{\la}(\la)]_t$, proving the Theorem for the case $\kappa(\la)=0$.

\medskip
Now we suppose that $\kappa(\bs{\lambda})=2$. By Lemma \ref{cantidad M_n} we get $|M_{n}(\bs{\lambda})|=6$. In this setting, we have the following three possibilities for the order $\succeq$ on $M_{n}(\bs{\lambda})$
\begin{align*}
\bs{\lambda}_1 \prec \bs{\lambda}_2 \prec \bs{\lambda}_3\prec \bs{\lambda}_4 \prec \bs{\lambda}_5\prec \bs{\lambda}_6 \\
\bs{\lambda}_1 \prec \bs{\lambda}_3 \prec \bs{\lambda}_2\prec \bs{\lambda}_4 \prec \bs{\lambda}_5\prec \bs{\lambda}_6 \\
\bs{\lambda}_1 \prec \bs{\lambda}_2 \prec \bs{\lambda}_3\prec \bs{\lambda}_5 \prec \bs{\lambda}_4\prec \bs{\lambda}_6
\end{align*}
In this three cases, the theorem follows by a case to case analysis, we only prove the lemma for the (most interesting) last case. Figure below shows an example of this case for $\la=(0,16)$, $l=5$ and $m=2$. In the figure are also drawn all walks $w(\sub{s})$ with $\sub{s} \approx \sub{t}^{\la}$.

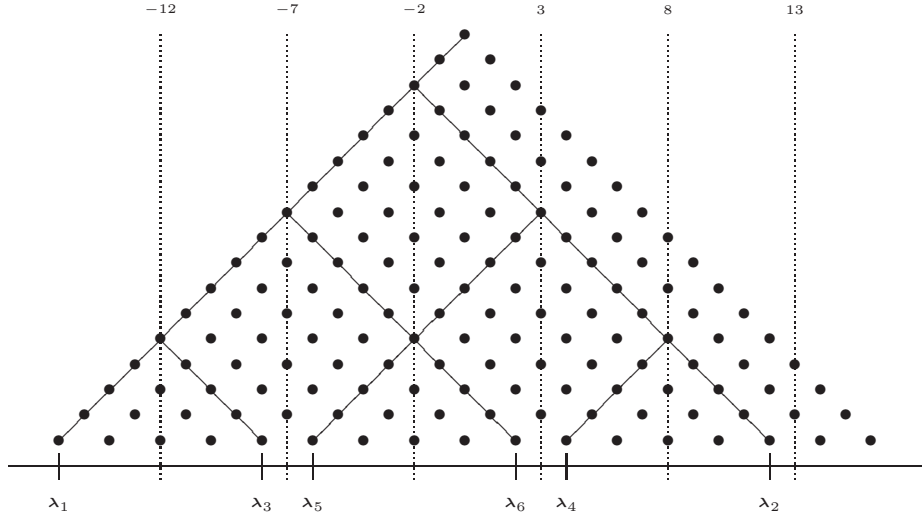
$$  \maximalesdos$$

By Corollary \ref{casi simple modules}, the modules $L_{\bs{\lambda}}(\bs{\lambda}_1)$ and $L_{\bs{\lambda}}(\bs{\lambda}_{5})$ are the unique (up degree shift) graded simple (non-zero) modules for $b_n(m,\la)$, furthermore by Remark \ref{dimension simples mayor que uno} these modules are one dimensional and pure of degree zero. Now,  Theorem \ref{gd standard} implies
$$   \begin{array}{lll}
     \dim_{t}\Delta_{\bs{\lambda}}(\bs{\lambda}_1)=1;   & \dim_{t}\Delta_{\bs{\lambda}}(\bs{\lambda}_2)=t; & \dim_{t}\Delta_{\bs{\lambda}}(\bs{\lambda}_3)=t;  \\
      \dim_{t}\Delta_{\bs{\lambda}}(\bs{\lambda}_4)=t^{2};  &  \dim_{t}\Delta_{\bs{\lambda}}(\bs{\lambda}_5)=t^2+1; & \dim_{t}\Delta_{\bs{\lambda}}(\bs{\lambda}_6)=t^3+t.
     \end{array}
 $$
By Theorem \ref{gdn in a cellular algebra} $[\Delta_{\la}(\la):L_{\la}(\la)]_t=1$. As in the previous case we analyze equation (\ref{relacion dimension cell y simples}) for the different values of $k$. For $k=2,3$ equation (\ref{relacion dimension cell y simples}) becomes
\begin{equation*}
t=  [\Delta_{\la}(\la_{2}):L_{\la}(\la)]_t    \quad \text{ and } \quad   t=  [\Delta_{\la}(\la_{3}):L_{\la}(\la)]_t
\end{equation*}
Next, if $k=4$ then equation (\ref{relacion dimension cell y simples}) becomes
\begin{equation*}
t^{2}=  [\Delta_{\la}(\la_{4}):L_{\la}(\la_5)]_t+  [\Delta_{\la}(\la_{4}):L_{\la}(\la)]_t
\end{equation*}
but it is straightforward to check that $\bs{\lambda}_{4} \not \in M_n(\bs{\lambda}_{5})$, hence by Lemma \ref{gdn zero}  $[\Delta_{\la}(\la_{4}):L_{\la}(\la_5)]_t=0$, thus $t^{2}=  [\Delta_{\la}(\la_{4}):L_{\la}(\la_5)]_t$. Now, for $k=5$ we have

\begin{equation*}
t^2+1=  [\Delta_{\la}(\la_{5}):L_{\la}(\la_5)]_t  +   [\Delta_{\la}(\la_{5}):L_{\la}(\la)]_t
\end{equation*}
but by Theorem \ref{gdn in a cellular algebra}(c) we have  $[\Delta_{\la}(\la_{5}):L_{\la}(\la_5)]_t=1$, so   $[\Delta_{\la}(\la_{5}):L_{\la}(\la_5)]_t=t^{2}$. Finally, for $k=6$ equation (\ref{relacion dimension cell y simples}) becomes
\begin{equation*}
t^3+t=  [\Delta_{\la}(\la_{6}):L_{\la}(\la_5)]_t  +   [\Delta_{\la}(\la_{6}):L_{\la}(\la)]_t
\end{equation*}
It is not hard to note that $\kappa(\bs{\lambda}_5)=0$ and  that  $M_n(\bs{\lambda}_5)=\{ \bs{\lambda}_5, \bs{\lambda}_6\}$, so we know by the first case analyzed in this proof that $[\Delta_{\bs{\lambda}}(\bs{\lambda}_6) : L_{\bs{\lambda}}(\bs{\lambda}_5)]_t=t$. Consequently, $[\Delta_{\la}(\la_{6}):L_{\la}(\la)]_t=t^{3}$. This completes the proof of the Theorem for the case $\kappa(\la)=2$.

\medskip
Now we can assume that $\kappa(\bs{\lambda}) \neq 0,2$. By Remark \ref{dimension simples mayor que uno}, $ L_{\bs{\lambda}}(\bs{\lambda})$  is the unique (up degree shift) one-dimensional graded simple module for $b_{n}(m,\bs{\lambda})$. Recall that $b_{n}(m,\la)$ is a positively graded cellular algebra, so we can consider the grading filtration for $\Delta_{\la}(\la_k)$, $\la_k \in M_n(\la)$. Now, by Theorem \ref{gd standard} and  Remark \ref{dimension simples mayor que uno} $\dim_{t}\Delta_{\la}(\la_k) \in \Z[t]$ is a monic polynomial with the non-leader coefficients greater than $1$. Thus, in the grading filtration of $\Delta_{\la}(\la_k)$ there is a unique quotient of dimension one. This quotient is pure of degree $\deg (\dim_{t}\Delta_{\la}(\la_k))$ (where here $\deg$ denotes the polynomial degree) and must be isomorphic (in the ungraded setting) to $L_{\la}(\la)$. Since $L_{\la}(\la)$ is pure of degree zero, if the grading filtration for $\Delta_{\la}(\la_k)$ is a graded composition series we have
 \begin{equation} \label{casi teorema main}
[\Delta_{\la}(\la_{k}):L_{\la}(\la)]_t = t^{\deg (\dim_{t}\Delta_{\la}(\la_k))}
\end{equation}
If the grading filtration for $\Delta_{\la}(\la_k)$ is not a graded composition series we can always add graded $b_n(m,\la)$-submodules of $\Delta_{\la}(\la_k)$ to the grading filtration in order to obtain a graded composition series. In a graded composition series obtained in this way we can also have only one graded composition factor of dimension one. Otherwise we obtain via Theorem \ref{gdn iguales} that the blob algebra $b_n(m)$ has a (ungraded) decomposition number greater than one. Therefore, (\ref{casi teorema main}) is still valid even if the grading filtration is not a graded composition series. Finally, from Theorem \ref{gd standard} we know that
\begin{equation}  \label{ultima}
\deg (\dim_{t}( \Delta_{\bs{\lambda}}(\bs{\lambda}_k)))= \left\{
                                          \begin{array}{ll}
                                            2j,   &  \text{ if } k=4j+1 ;    \\
                                            2j+1, &  \text{ if } k=4j+2;  \\
                                            2j+1, & \text{ if } k=4j+3;   \\
                                            2j+2, & \text{ if } k=4j+4.
                                          \end{array}
                                        \right.
 \end{equation}
Now the Theorem follows by combining (\ref{casi teorema main}) and (\ref{ultima}).
\end{dem}

\subsection{The wall case.}
In the previous subsection we determine the graded decomposition numbers for $\blob$, $[\Delta(\bs{\mu}): L(\bs{\lambda})]_t$, when $\la$ is not on a wall. Now we determine these numbers in the \emph{wall case}. So in this subsection we assume that $\bs{\lambda}$ is on a wall. For brevity, the results shall be presented without proof since the same series of arguments used to the non-wall case working here with minimal changes. If $\bs{\lambda}$ is on one of the walls of the fundamental alcove we have $M_{n}(\bs{\lambda})=\{ \bs{\lambda}\}$ so by Theorem \ref{gdn in a cellular algebra} and Lemma \ref{gdn zero} we obtain $[\Delta(\bs{\mu}): L(\bs{\lambda})]_t=\delta_{\bs{\mu \lambda}}$, where $\delta_{\bs{\mu \lambda}}$ is a Kronecker delta. So we can assume that $\bs{\lambda}$ is not a wall of the fundamental alcove. As in the non-wall case, define  $\kappa(\bs{\lambda})$ as the number of alcoves between $\la$ and $0$. Similarly, we index $M_n(\bs{\lambda})$ and define two-column partitions of $\kappa(\la)$, $ \mu_{i}^{j}$, by the same rules used for the non-wall case.
The following theorems correspond to Theorem \ref{gd standard} and Theorem \ref{main theorem} in the \emph{non-wall} case.
\begin{teo}  \label{gd standard wall}
Let $\bs{\lambda} \in \bi{n}$. The graded dimension of $\Delta_{\bs{\lambda}}(\bs{\lambda}_k)$, for $\la_k \in M_{n}(\la) $, is completely determined by the formulas:
$$ \begin{array}{l}
 (a)  \qquad   \dim_t(\Delta_{\bs{\lambda}}(\bs{\lambda}_{2j+1}) )= \sum_{i=0}^{j} c_{i}t^{2i} \\
   \\
(b) \qquad \dim_t(\Delta_{\bs{\lambda}}(\bs{\lambda}_{2j+2}) )= t \dim_t(\Delta_{\bs{\lambda}}(\bs{\lambda}_{2j+1}) )   \end{array}
$$
where $c_i=|\Std(\mu_{i}^{j})|$.
\end{teo}

\begin{teo} Let $\bs{\lambda} \in \bi{n}$. For $\la_k\in M_{n}(\la)$ we have
 \begin{equation*} \label{gdn answer wall}
 [\Delta(\bs{\lambda}_k) : L(\bs{\lambda})]_t= t^{k-1}
 \end{equation*}
\end{teo}

\begin{scriptsize}

\sc Instituto de Matem\'atica y F\'isica, Universidad de Talca, Chile, mail dplaza@inst-mat.utalca.cl
\end{scriptsize}

\end{document}